\documentstyle[amstex,amssymb,11pt]{article}
  \input xypic
 \newcommand{\lon}{\longrightarrow}
 
 \newcommand{\rar}{\rightarrow}
 
 \newcommand{\Proof}{{\bf Proof}.\, }

 \newcommand{\Z}{{\Bbb Z}}
 
 \newcommand{\p}{{\partial}}

 \newcommand{\R}{{\Bbb R}}
 \newcommand{\ot}{\otimes}

 \newcommand{\Id}{\mbox{Id}}
\newcommand{\LieBi}{{\mathsf Lie^1Bi }}
 \newcommand{\Beq}{\begin{equation}}
 \newcommand{\Eeq}{\end{equation}}
 \newcommand{\Beqr}{\begin{eqnarray}}
 \newcommand{\Eeqr}{\end{eqnarray}}
 \newcommand{\Beqrn}{\begin{eqnarray*}}
 \newcommand{\Eeqrn}{\end{eqnarray*}}
 \newcommand{\Ba}{\begin{array}}
 \newcommand{\Ea}{\end{array}}
 \newcommand{\Bi}{\begin{itemize}}
 \newcommand{\Ei}{\end{itemize}}
 \newcommand{\Bc}{\begin{center}}
 \newcommand{\Ec}{\end{center}}

 \newcommand{\f}{{\cal O}}

 \newcommand{\cE}{{\cal E}}

 \newcommand{\caL}{{\cal L}}
 \newcommand{\cM}{{\cal M}}
 \newcommand{\cP}{{\cal P}}

 \newcommand{\al}{\alpha}
 \newcommand{\be}{\beta}

 \newcommand{\Ga}{\Gamma}

 \newcommand{\om}{\omega}

 \newcommand{\Img}{{\mathsf Im}\, }

 \newcommand{\sip}{\smallskip}
 \newcommand{\bip}{\bigskip}

\begin{document}

 \sloppy

 \title{PROP profile of Poisson geometry}
 \author{ S.A.\ Merkulov}
 \date{}
 \maketitle

\hspace{77mm}{\small ``The genetic code appears to be universal; \ldots''

\hspace{120mm} {\small\em Britannica.}

\bip

{\bf 0. Introduction.} The first instances of algebraic and
topological strongly homotopy, or infinity, structures have been
discovered by Stasheff \cite{St} long ago. Since that time
infinities have acquired a prominent role in algebraic
topology and homological algebra. 
We argue in this paper that some classical local geometries are of
infinity origin, i.e.\ their smooth formal germs are (homotopy)
representations of cofibrant PROPs
$\cP_\infty$ in spaces
concentrated in degree zero; in particular, they admit natural
infinity generalizations when one considers homotopy
representations of
$\cP_\infty$ in generic differential graded (dg) spaces. The simplest
manifestation of this phenomenon is provided by the Poisson geometry
(or even by smooth germs of tensor fields!) and is the main theme of the present paper.
Another example is discussed in \cite{Me2}.
The
PROPs $\cP_\infty$ are minimal resolutions of PROPs $\cP$ which are 
graph spaces built from very few basic elements, {\em genes}, subject to
simple engineering rules. Thus to a local geometric structure one
can associate a kind of a code, {\em genome}, which specifies it
uniquely and opens a new window of opportunities of attacking  differential
 geometric problems
with the powerful tools of homological algebra. In particular, the genetic
code of Poisson geometry discovered in this paper
 has been used in
 \cite{Me3} to give a new short proof of 
Kontsevich's \cite{Ko} deformation quantization theorem.

\sip

Formal germs of geometric structures discussed in this
paper are {\em pointed} \, in the sense that they vanish at the
distinguished point.
This is the usual price one pays for working
with (di)operads without ``zero terms" (as is often done in the
literature).  As structural equations behind the particular geometries
we study in this paper
are homogeneous, this restriction poses no problem: say, a generic non-pointed
Poisson structure, $\nu$, in $\R^n$ can be identified with the pointed
one, $\hbar \nu$, in  $\R^{n+1}$, $\hbar$ being the extra coordinate.

\sip

We introduce in this paper a
dg free dioperad whose generic representations in  a graded vector
space $V$ can  be identified with  pointed solutions of the Maurer-Cartan equations
in the Lie algebra of polyvector fields on the formal manifold associated
with $V$. The cohomology of this dioperad can not be computed directly. Instead
one has to rely on some fine mathematics such as Koszulness \cite{GK,G}
 and distributive laws \cite{Mar,G}. One of the main results of this paper is a
proof of Theorem 3.2 which identifies the cohomology of that dg free dioperad
with a surprisingly small dioperad, $\LieBi$, of Lie 1-bialgebras,   
which are almost identical to the dioperad, ${\mathsf LieBi}$,
 of usual Lie bialgebras except
that degree of generating Lie and coLie operations differ by 1
(compare with Gerstenhaber versus Poisson algebras).
The dioperad $\LieBi$ is proven to be Koszul. We use the resulting geometric
interpretation of $\LieBi_\infty$ algebras to give their homotopy classification
(see Theorem 3.4.5) which is an extension of Kontsevich's homotopy classification
\cite{Ko}
of ${\mathsf L}_\infty$ algebras.

\sip

As a side remark we also discuss graph and geometric interpretations
of strongly homotopy
Lie bialgebras using Koszulness of the
latter which was established in \cite{G}.


\bip

 {\bf 1.~{ Geometry} $\Rightarrow$ PROP profile  $\Rightarrow$ {Geometry}$_\infty$.}
Let $\cP$ be an operad, or a dioperad, or even a PROP admitting a
minimal dg resolution. Let $\cP{\mathsf Alg}$ be
the category of finite dimensional dg $\cP$-algebras, and  ${\mathbf D}(\cP{\mathsf Alg})$
the associated
derived category (which we understand here as the homotopy category of $\cP_{\infty}$-algebras,
$\cP_\infty$ being the minimal resolution of $\cP$).

\sip

For any locally defined
geometric structure ${\sf Geom}$ (say, Poisson, Riemann, K\"ahler, etc.) it makes sense
talking  about the category of formal  ${\sf Geom}$-manifolds. Its objects
are  formal pointed manifolds (non-canonically isomorphic to   $(\R^n, 0)$ for some  $n$)
together with a germ
of formal ${\sf Geom}$-structure
at  the distinguished point.

\bip

{\bf 1.1. Definition.} The operad/dioperad/PROP
$\cP$ is called a PROP-{\em profile},
or {\em genome},  of
a geometric structure ${\sf Geom}$ if
\Bi
\item
 the category
of formal  ${\sf Geom}$-manifolds is equivalent to a full subcategory of the derived category
 ${\mathbf D}(\cP{\mathsf Alg})$ , and
\item there is no sub-(di)operad of $\cP$ having the above property.
\Ei


\bip

{\bf 1.2. Definition.} If $\cP$ is a PROP-profile of a geometric structure
  ${\sf Geom}$, then a generic object
of  ${\mathbf D}(\cP{\mathsf Alg})$  is called a formal ${\sf Geom}_\infty$-manifold.

\bip

 Presumably, ${\sf Geom}_\infty$-structure is what one gets from  ${\sf Geom}$
by means of the extended deformation theory.

\sip

Local geometric structures are often  non-trivial and  complicated
creatures --- the general solution of the associated defining
system of nonlinear differential equations is not available; it is
often a very hard job just to show  existence of non-trivial
solutions. Nevertheless, if such a structure  ${\sf Geom}$ admits
a PROP-profile, $\cP= Free(\cE)/Ideal$\footnote{Any
operad/dioperad/etc. can be represented as a quotient of the free
operad/dioperad/etc., $Free(\cE)$ generated by a collection of
$\Sigma_m$-left/$\Sigma_n$-right modules $\cE=\{\cE(m,n)\}_{m,n\geq 1}$, by an $Ideal$.
Often there exists a canonical, ``common factors canceled out'',
representation like this.}, then  ${\sf Geom}$  can be
non-ambiguously  characterized by its ``genetic code'': {\em
genes}\, are, by definition, the generators of $\cE$, and the {\em
engineering rules}\, are, by definition, the generators of
$Ideal$. And that code can be surprisingly simple, as examples 1.3-1.5 and
the table below illustrate.

\bip

\begin{table}
\begin{center}
\begin{tabular}{|c|c|c|}
\multicolumn{3}{c}{}\vspace{2 mm}\\
\hline
              & generic representation     & generic representation \\
 Genome $\cP$ & of $\cP_\infty$ in $\R^n$  &  of $\cP_\infty$ in a graded\\
              &                            &  vector space ${V}$          \\
\hline
\hline

\begin{tabular}{c} \vspace{- 3mm}\\
                     $\cP$ is the  $G$-operad\vspace{1 mm}\\
                      Genes:

$
 \begin{xy}
 <0mm,0.66mm>*{};<0mm,3mm>*{}**@{-},
 <0.39mm,-0.39mm>*{};<2.2mm,-2.2mm>*{}**@{-},
 <-0.35mm,-0.35mm>*{};<-2.2mm,-2.2mm>*{}**@{-},
 <0mm,0mm>*{\circ};<0mm,0mm>*{}**@{},
 \end{xy}$\  , \,
$
 \begin{xy}
 <0mm,0.66mm>*{};<0mm,3mm>*{}**@{-},
 <0.39mm,-0.39mm>*{};<2.2mm,-2.2mm>*{}**@{-},
 <-0.35mm,-0.35mm>*{};<-2.2mm,-2.2mm>*{}**@{-},
 <0mm,0mm>*{\bullet};<0mm,0mm>*{}**@{},
 \end{xy}$
 \vspace{1 mm}\\
                      Engineering rules: \vspace{1 mm}
$\begin{xy}
 <0mm,0mm>*{\circ};<0mm,0mm>*{}**@{},
 <0mm,0.69mm>*{};<0mm,3.0mm>*{}**@{-},
 <0.39mm,-0.39mm>*{};<1.9mm,-1.9mm>*{}**@{-},
 <-0.35mm,-0.35mm>*{};<-2.2mm,-2.2mm>*{}**@{-},
 <2.4mm,-2.4mm>*{\circ};<2.4mm,-2.4mm>*{}**@{},
 <2.0mm,-2.8mm>*{};<0mm,-4.9mm>*{}**@{-},
 <2.8mm,-2.9mm>*{};<4.6mm,-4.9mm>*{}**@{-},
 \end{xy}$
-
$
 \begin{xy}
 <0mm,0mm>*{\circ};<0mm,0mm>*{}**@{},
 <0mm,0.69mm>*{};<0mm,3.0mm>*{}**@{-},
 <0.39mm,-0.39mm>*{};<2.2mm,-2.2mm>*{}**@{-},
 <-0.35mm,-0.35mm>*{};<-1.9mm,-1.9mm>*{}**@{-},
 <-2.4mm,-2.4mm>*{\circ};<-2.4mm,-2.4mm>*{}**@{},
 <-2.0mm,-2.8mm>*{};<0mm,-4.9mm>*{}**@{-},
 <-2.8mm,-2.9mm>*{};<-4.7mm,-4.9mm>*{}**@{-},
 \end{xy}
$
$=0$ \\
$$
 \begin{xy}
 <0mm,0mm>*{\bullet};<0mm,0mm>*{}**@{},
 <0mm,0.69mm>*{};<0mm,3.0mm>*{}**@{-},
 <0.39mm,-0.39mm>*{};<1.9mm,-1.9mm>*{}**@{-},
 <-0.35mm,-0.35mm>*{};<-4.2mm,-4.9mm>*{}**@{-},
 <2.4mm,-2.4mm>*{\bullet};<2.4mm,-2.4mm>*{}**@{},
 <2.0mm,-2.8mm>*{};<0mm,-4.9mm>*{}**@{-},
 <2.8mm,-2.9mm>*{};<4.6mm,-4.9mm>*{}**@{-},
 \end{xy}
+
\,
\begin{xy}
 <0mm,0mm>*{\bullet};<0mm,0mm>*{}**@{},
 <0mm,0.69mm>*{};<0mm,3.0mm>*{}**@{-},
 <0.39mm,-0.39mm>*{};<1.9mm,-1.9mm>*{}**@{-},
 <0mm,0mm>*{};<-1.5mm,-1.5mm>*{}**@{-},
 <-1.5mm,-1.5mm>*{};<2.4mm,-4.9mm>*{}**@{-},
 <2.4mm,-2.4mm>*{\bullet};<2.4mm,-2.4mm>*{}**@{},
 <2.0mm,-2.8mm>*{};<0mm,-4.9mm>*{}**@{-},
 <2.8mm,-2.9mm>*{};<4.8mm,-4.9mm>*{}**@{-},
 \end{xy}
+
 \begin{xy}
 <0mm,0mm>*{\bullet};<0mm,0mm>*{}**@{},
 <0mm,0.69mm>*{};<0mm,3.0mm>*{}**@{-},
 <0.39mm,-0.39mm>*{};<4.3mm,-4.9mm>*{}**@{-},
 <-0.35mm,-0.35mm>*{};<-1.9mm,-1.9mm>*{}**@{-},
 <-2.4mm,-2.4mm>*{\bullet};<-2.4mm,-2.4mm>*{}**@{},
 <-2.0mm,-2.8mm>*{};<0mm,-4.9mm>*{}**@{-},
 <-2.8mm,-2.9mm>*{};<-4.7mm,-4.9mm>*{}**@{-},
 \end{xy}
=\ 0
$$
\vspace{1mm}\\
$$
\begin{xy}
 <0mm,0mm>*{\bullet};<0mm,0mm>*{}**@{},
 <0mm,0.69mm>*{};<0mm,3.0mm>*{}**@{-},
 <0.39mm,-0.39mm>*{};<1.9mm,-1.9mm>*{}**@{-},
 <-0.35mm,-0.35mm>*{};<-4.2mm,-4.9mm>*{}**@{-},
 <2.4mm,-2.4mm>*{\circ};<2.4mm,-2.4mm>*{}**@{},
 <2.0mm,-2.8mm>*{};<0mm,-4.9mm>*{}**@{-},
 <2.8mm,-2.9mm>*{};<4.6mm,-4.9mm>*{}**@{-},
 \end{xy}
-
\,
\begin{xy}
 <0mm,0mm>*{\circ};<0mm,0mm>*{}**@{},
 <0mm,0.69mm>*{};<0mm,3.0mm>*{}**@{-},
 <0.39mm,-0.39mm>*{};<1.9mm,-1.9mm>*{}**@{-},
 <-0.39mm,-0.39mm>*{};<-1.7mm,-1.7mm>*{}**@{-},
 <-1.7mm,-1.7mm>*{};<2.4mm,-4.9mm>*{}**@{-},
 <2.4mm,-2.4mm>*{\bullet};<2.4mm,-2.4mm>*{}**@{},
 <2.0mm,-2.8mm>*{};<0mm,-4.9mm>*{}**@{-},
 <2.8mm,-2.9mm>*{};<4.8mm,-4.9mm>*{}**@{-},
 \end{xy}
-
 \begin{xy}
 <0mm,0mm>*{\circ};<0mm,0mm>*{}**@{},
 <0mm,0.69mm>*{};<0mm,3.0mm>*{}**@{-},
 <0.39mm,-0.39mm>*{};<4.3mm,-4.9mm>*{}**@{-},
 <-0.35mm,-0.35mm>*{};<-1.9mm,-1.9mm>*{}**@{-},
 <-2.4mm,-2.4mm>*{\bullet};<-2.4mm,-2.4mm>*{}**@{},
 <-2.0mm,-2.8mm>*{};<0mm,-4.9mm>*{}**@{-},
 <-2.8mm,-2.9mm>*{};<-4.7mm,-4.9mm>*{}**@{-},
 \end{xy}
=\ 0
$$
\vspace{2mm} \\

                       \end{tabular}

  & \begin{tabular}{c} \vspace{- 3mm}\\
                      smooth formal   \\  Hertling-Manin\\
                      structure in $\R^n$ \cite{HM}\vspace{1 mm}\\
                       \vspace{1 mm}
                       \end{tabular}
&
\begin{tabular}{c} \vspace{- 3mm}\\
                       smooth formal  \\ Hertling-Manin$_{\infty}$\\
                      structure in $\hat{V}$ \cite{Me}\vspace{1 mm}\\
                       \vspace{1 mm}
                       \end{tabular}
 \\
\hline
\hline
\begin{tabular}{c} \vspace{- 3mm}\\
                   $\cP$ is the dioperad ${\mathsf TF}$\vspace{1 mm}\\
                      Genes:

$
 \begin{xy}
 <0mm,-0.55mm>*{};<0mm,-2.5mm>*{}**@{-},
 <0.5mm,0.5mm>*{};<2.2mm,2.2mm>*{}**@{-},
 <-0.48mm,0.48mm>*{};<-2.2mm,2.2mm>*{}**@{-},
 <0mm,0mm>*{\circ};<0mm,0mm>*{}**@{},
 \end{xy}
$ , \ $
 \begin{xy}
 <0mm,0.66mm>*{};<0mm,3mm>*{}**@{-},
 <0.39mm,-0.39mm>*{};<2.2mm,-2.2mm>*{}**@{-},
 <-0.35mm,-0.35mm>*{};<-2.2mm,-2.2mm>*{}**@{-},
 <0mm,0mm>*{\bullet};<0mm,0mm>*{}**@{},
 \end{xy}$
 \vspace{1 mm}\\
                      Rules:  \ \ \ \ \ \ \ \ \ \ \ \vspace{1 mm}\\
$$
 \begin{xy}
 <0mm,0mm>*{\bullet};<0mm,0mm>*{}**@{},
 <0mm,0.69mm>*{};<0mm,3.0mm>*{}**@{-},
 <0.39mm,-0.39mm>*{};<1.9mm,-1.9mm>*{}**@{-},
 <-0.35mm,-0.35mm>*{};<-4.2mm,-4.9mm>*{}**@{-},
 <2.4mm,-2.4mm>*{\bullet};<2.4mm,-2.4mm>*{}**@{},
 <2.0mm,-2.8mm>*{};<0mm,-4.9mm>*{}**@{-},
 <2.8mm,-2.9mm>*{};<4.6mm,-4.9mm>*{}**@{-},
 \end{xy}
+ \,
\begin{xy}
 <0mm,0mm>*{\bullet};<0mm,0mm>*{}**@{},
 <0mm,0.69mm>*{};<0mm,3.0mm>*{}**@{-},
 <0.39mm,-0.39mm>*{};<1.9mm,-1.9mm>*{}**@{-},
 <0mm,0mm>*{};<-1.5mm,-1.5mm>*{}**@{-},
 <-1.5mm,-1.5mm>*{};<2.4mm,-4.9mm>*{}**@{-},
 <2.4mm,-2.4mm>*{\bullet};<2.4mm,-2.4mm>*{}**@{},
 <2.0mm,-2.8mm>*{};<0mm,-4.9mm>*{}**@{-},
 <2.8mm,-2.9mm>*{};<4.8mm,-4.9mm>*{}**@{-},
 \end{xy}
+
 \begin{xy}
 <0mm,0mm>*{\bullet};<0mm,0mm>*{}**@{},
 <0mm,0.69mm>*{};<0mm,3.0mm>*{}**@{-},
 <0.39mm,-0.39mm>*{};<4.3mm,-4.9mm>*{}**@{-},
 <-0.35mm,-0.35mm>*{};<-1.9mm,-1.9mm>*{}**@{-},
 <-2.4mm,-2.4mm>*{\bullet};<-2.4mm,-2.4mm>*{}**@{},
 <-2.0mm,-2.8mm>*{};<0mm,-4.9mm>*{}**@{-},
 <-2.8mm,-2.9mm>*{};<-4.7mm,-4.9mm>*{}**@{-},
 \end{xy}
=\ 0
$$
\vspace{2mm}\\
$$
 \begin{xy}
 <0mm,2.47mm>*{};<0mm,-0.5mm>*{}**@{-},
 <0.5mm,3.5mm>*{};<2.2mm,5.2mm>*{}**@{-},
 <-0.48mm,3.48mm>*{};<-2.2mm,5.2mm>*{}**@{-},
 <0mm,3mm>*{\circ};<0mm,3mm>*{}**@{},
  <0mm,-0.8mm>*{\bullet};<0mm,-0.8mm>*{}**@{},
<0mm,-0.8mm>*{};<-2.2mm,-3.5mm>*{}**@{-},
 <0mm,-0.8mm>*{};<2.2mm,-3.5mm>*{}**@{-},
\end{xy}
 =
\begin{xy}
 <0mm,-1.3mm>*{};<0mm,-3.5mm>*{}**@{-},
 <0.38mm,-0.2mm>*{};<2.2mm,2.2mm>*{}**@{-},
 <-0.38mm,-0.2mm>*{};<-2.2mm,2.2mm>*{}**@{-},
<0mm,-0.8mm>*{\circ};-<0mm,0.8mm>*{}**@{},
 <-2.25mm,2.2mm>*{};<-2.2mm,5.2mm>*{}**@{-},
 <2.4mm,2.4mm>*{\bullet};<2.4mm,2.4mm>*{}**@{},
 <2.5mm,2.3mm>*{};<4.4mm,-0.8mm>*{}**@{-},
 <4.4mm,-0.8mm>*{};<4.4mm,-3.5mm>*{}**@{-},
 <2.4mm,2.5mm>*{};<2.4mm,5.2mm>*{}**@{-},
 \end{xy}
 +
\begin{xy}
<3mm,-1.3mm>*{};<3mm,-3.5mm>*{}**@{-},
<3.4mm,-0.2mm>*{};<5.2mm,2.2mm>*{}**@{-},
<2.6mm,-0.2mm>*{};<0.4mm,2.2mm>*{}**@{-},
<3mm,-0.8mm>*{\circ};<3mm,-0.8mm>*{}**@{},
<5.2mm,2.2mm>*{};<5.2mm,5.2mm>*{}**@{-},
<0.3mm,2.4mm>*{\bullet};<0.3mm,2.4mm>*{}**@{},
<0.1mm,2.3mm>*{};<-2.2mm,-0.8mm>*{}**@{-},
<-2.2mm,-0.8mm>*{};<-2.2mm,-3.5mm>*{}**@{-},
<0.3mm,2.5mm>*{};<0.3mm,5.2mm>*{}**@{-},
\end{xy}
 +
\begin{xy}
<3mm,-0.8mm>*{\circ};<3mm,-0.8mm>*{}**@{},
<3mm,-1.3mm>*{};<3mm,-3.5mm>*{}**@{-},
<3.7mm,-0.5mm>*{};<5.7mm,1.2mm>*{}**@{-},
<5.7mm,1.2mm>*{};<3.2mm,2.4mm>*{}**@{-},
<2.5mm,-0.5mm>*{};<-0.5mm,2.4mm>*{}**@{-},
<-0.5mm,2.4mm>*{};<-0.5mm,5.2mm>*{}**@{-},
<3mm,2.4mm>*{\bullet};<3mm,2.4mm>*{}**@{},
<3.3mm,2.9mm>*{};<1.4mm,1.0mm>*{}**@{-},
<0.9mm,0.5mm>*{};<-0.5mm,-0.8mm>*{}**@{-},
<-0.5mm,-0.8mm>*{};<-0.5mm,-3.5mm>*{}**@{-},
<3mm,2.4mm>*{};<3mm,5.2mm>*{}**@{-},
\end{xy}
 +
\begin{xy}
<3mm,-0.8mm>*{\circ};<3mm,-0.8mm>*{}**@{},
<3mm,-1.3mm>*{};<3mm,-3.5mm>*{}**@{-},
<2.3mm,-0.5mm>*{};<0.3mm,1.2mm>*{}**@{-},
<0.3mm,1.2mm>*{};<2.4mm,2.4mm>*{}**@{-},
<3.5mm,-0.5mm>*{};<7mm,2.4mm>*{}**@{-},
<7mm,2.4mm>*{};<7mm,5.2mm>*{}**@{-},
<3mm,2.4mm>*{\bullet};<3mm,2.4mm>*{}**@{},
<2.7mm,2.9mm>*{};<4.8mm,1.0mm>*{}**@{-},
<5.4mm,0.5mm>*{};<7mm,-0.8mm>*{}**@{-},
<7mm,-0.8mm>*{};<7mm,-3.5mm>*{}**@{-},
<3mm,2.4mm>*{};<3mm,5.2mm>*{}**@{-},
\end{xy}
$$

\vspace{1mm}

                       \end{tabular}

  & \begin{tabular}{l} \vspace{1mm}\\
                      smooth formal \\ section
                      of $\ot^2T_{\R^n}$ \\
                      (variants: of $\wedge^2T_{\R^n}$\\
                      or of  $\odot^2T_{\R^n}$) \\
                      vanishing at $0$  \vspace{1 mm}\\
                       \vspace{1 mm}

                      \end{tabular}
&
\begin{tabular}{l} \vspace{- 3mm}\\
                      structure, $(\hat{V},\eth\in T_{\hat{V}}$),\\
                      of a smooth dg manifold \\
                       together with a
                      smooth  \\ section $\phi$
                      of $\ot^2T_{\hat{V}}$ \\
                       (variants: of $\wedge^2T_{\hat{V}}$\\
                      or of  $\odot^2T_{\hat{V}}$) \\
                      vanishing at $0$ and\\
                      satisfying $Lie_\eth \phi=0$.
                       \\
                       \vspace{1 mm}\\
                       \end{tabular}

 \\
\hline
\hline
\begin{tabular}{c} \vspace{- 3mm}\\
                   $\cP$ is the dioperad $\LieBi$\vspace{1 mm}\\
                      Genes:

$
 \begin{xy}
 <0mm,-0.55mm>*{};<0mm,-2.5mm>*{}**@{-},
 <0.5mm,0.5mm>*{};<2.2mm,2.2mm>*{}**@{-},
 <-0.48mm,0.48mm>*{};<-2.2mm,2.2mm>*{}**@{-},
 <0mm,0mm>*{\circ};<0mm,0mm>*{}**@{},
 \end{xy}
$ , \
$
 \begin{xy}
 <0mm,0.66mm>*{};<0mm,3mm>*{}**@{-},
 <0.39mm,-0.39mm>*{};<2.2mm,-2.2mm>*{}**@{-},
 <-0.35mm,-0.35mm>*{};<-2.2mm,-2.2mm>*{}**@{-},
 <0mm,0mm>*{\bullet};<0mm,0mm>*{}**@{},
 \end{xy}$
 \vspace{1 mm}\\
                      Rules: \vspace{1 mm}
%
$$
\begin{xy}
 <0mm,0mm>*{\circ};<0mm,0mm>*{}**@{},
 <0mm,-0.49mm>*{};<0mm,-3.0mm>*{}**@{-},
 <0.49mm,0.49mm>*{};<3.9mm,4.9mm>*{}**@{-},
 <-0.5mm,0.5mm>*{};<-1.9mm,1.9mm>*{}**@{-},
 <-2.3mm,2.3mm>*{\circ};<-2.3mm,2.3mm>*{}**@{},
 <-1.8mm,2.8mm>*{};<0mm,4.9mm>*{}**@{-},
 <-2.8mm,2.9mm>*{};<-4.6mm,4.9mm>*{}**@{-},
 \end{xy}
-
\begin{xy}
 <0mm,0mm>*{\circ};<0mm,0mm>*{}**@{},
 <0mm,-0.49mm>*{};<0mm,-3.0mm>*{}**@{-},
 <0.49mm,0.49mm>*{};<1.8mm,1.8mm>*{}**@{-},
 <1.8mm,1.8mm>*{};<-2.1mm,4.9mm>*{}**@{-},
 <-0.5mm,0.5mm>*{};<-1.9mm,1.9mm>*{}**@{-},
 <-2.3mm,2.3mm>*{\circ};<-2.3mm,2.3mm>*{}**@{},
 <-1.8mm,2.8mm>*{};<0mm,4.9mm>*{}**@{-},
 <-2.8mm,2.9mm>*{};<-4.6mm,4.9mm>*{}**@{-},
 \end{xy}
-
\begin{xy}
 <0mm,0mm>*{\circ};<0mm,0mm>*{}**@{},
 <0mm,-0.49mm>*{};<0mm,-3.0mm>*{}**@{-},
 <0.49mm,0.49mm>*{};<1.9mm,1.9mm>*{}**@{-},
 <-0.4mm,0.4mm>*{};<-3.2mm,4.9mm>*{}**@{-},
 <2.3mm,2.3mm>*{\circ};<2.3mm,2.3mm>*{}**@{},
 <1.9mm,2.8mm>*{};<0.5mm,4.9mm>*{}**@{-},
 <2.8mm,2.9mm>*{};<4.6mm,4.9mm>*{}**@{-},
 \end{xy}
=\, 0
$$
 \vspace{1mm} \\
$$
 \begin{xy}
 <0mm,0mm>*{\bullet};<0mm,0mm>*{}**@{},
 <0mm,0.69mm>*{};<0mm,3.0mm>*{}**@{-},
 <0.39mm,-0.39mm>*{};<1.9mm,-1.9mm>*{}**@{-},
 <-0.35mm,-0.35mm>*{};<-4.2mm,-4.9mm>*{}**@{-},
 <2.4mm,-2.4mm>*{\bullet};<2.4mm,-2.4mm>*{}**@{},
 <2.0mm,-2.8mm>*{};<0mm,-4.9mm>*{}**@{-},
 <2.8mm,-2.9mm>*{};<4.6mm,-4.9mm>*{}**@{-},
 \end{xy}
+
\,
\begin{xy}
 <0mm,0mm>*{\bullet};<0mm,0mm>*{}**@{},
 <0mm,0.69mm>*{};<0mm,3.0mm>*{}**@{-},
 <0.39mm,-0.39mm>*{};<1.9mm,-1.9mm>*{}**@{-},
 <0mm,0mm>*{};<-1.5mm,-1.5mm>*{}**@{-},
 <-1.5mm,-1.5mm>*{};<2.4mm,-4.9mm>*{}**@{-},
 <2.4mm,-2.4mm>*{\bullet};<2.4mm,-2.4mm>*{}**@{},
 <2.0mm,-2.8mm>*{};<0mm,-4.9mm>*{}**@{-},
 <2.8mm,-2.9mm>*{};<4.8mm,-4.9mm>*{}**@{-},
 \end{xy}
+
 \begin{xy}
 <0mm,0mm>*{\bullet};<0mm,0mm>*{}**@{},
 <0mm,0.69mm>*{};<0mm,3.0mm>*{}**@{-},
 <0.39mm,-0.39mm>*{};<4.3mm,-4.9mm>*{}**@{-},
 <-0.35mm,-0.35mm>*{};<-1.9mm,-1.9mm>*{}**@{-},
 <-2.4mm,-2.4mm>*{\bullet};<-2.4mm,-2.4mm>*{}**@{},
 <-2.0mm,-2.8mm>*{};<0mm,-4.9mm>*{}**@{-},
 <-2.8mm,-2.9mm>*{};<-4.7mm,-4.9mm>*{}**@{-},
 \end{xy}
=\ 0
$$
\vspace{2mm}\\
$$
 \begin{xy}
 <0mm,2.47mm>*{};<0mm,-0.5mm>*{}**@{-},
 <0.5mm,3.5mm>*{};<2.2mm,5.2mm>*{}**@{-},
 <-0.48mm,3.48mm>*{};<-2.2mm,5.2mm>*{}**@{-},
 <0mm,3mm>*{\circ};<0mm,3mm>*{}**@{},
  <0mm,-0.8mm>*{\bullet};<0mm,-0.8mm>*{}**@{},
<0mm,-0.8mm>*{};<-2.2mm,-3.5mm>*{}**@{-},
 <0mm,-0.8mm>*{};<2.2mm,-3.5mm>*{}**@{-},
\end{xy}
 =
\begin{xy}
 <0mm,-1.3mm>*{};<0mm,-3.5mm>*{}**@{-},
 <0.38mm,-0.2mm>*{};<2.2mm,2.2mm>*{}**@{-},
 <-0.38mm,-0.2mm>*{};<-2.2mm,2.2mm>*{}**@{-},
<0mm,-0.8mm>*{\circ};-<0mm,0.8mm>*{}**@{},
 <-2.25mm,2.2mm>*{};<-2.2mm,5.2mm>*{}**@{-},
 <2.4mm,2.4mm>*{\bullet};<2.4mm,2.4mm>*{}**@{},
 <2.5mm,2.3mm>*{};<4.4mm,-0.8mm>*{}**@{-},
 <4.4mm,-0.8mm>*{};<4.4mm,-3.5mm>*{}**@{-},
 <2.4mm,2.5mm>*{};<2.4mm,5.2mm>*{}**@{-},
 \end{xy}
 +
\begin{xy}
<3mm,-1.3mm>*{};<3mm,-3.5mm>*{}**@{-},
<3.4mm,-0.2mm>*{};<5.2mm,2.2mm>*{}**@{-},
<2.6mm,-0.2mm>*{};<0.4mm,2.2mm>*{}**@{-},
<3mm,-0.8mm>*{\circ};<3mm,-0.8mm>*{}**@{},
<5.2mm,2.2mm>*{};<5.2mm,5.2mm>*{}**@{-},
<0.3mm,2.4mm>*{\bullet};<0.3mm,2.4mm>*{}**@{},
<0.1mm,2.3mm>*{};<-2.2mm,-0.8mm>*{}**@{-},
<-2.2mm,-0.8mm>*{};<-2.2mm,-3.5mm>*{}**@{-},
<0.3mm,2.5mm>*{};<0.3mm,5.2mm>*{}**@{-},
\end{xy}
 +
\begin{xy}
<3mm,-0.8mm>*{\circ};<3mm,-0.8mm>*{}**@{},
<3mm,-1.3mm>*{};<3mm,-3.5mm>*{}**@{-},
<3.7mm,-0.5mm>*{};<5.7mm,1.2mm>*{}**@{-},
<5.7mm,1.2mm>*{};<3.2mm,2.4mm>*{}**@{-},
<2.5mm,-0.5mm>*{};<-0.5mm,2.4mm>*{}**@{-},
<-0.5mm,2.4mm>*{};<-0.5mm,5.2mm>*{}**@{-},
<3mm,2.4mm>*{\bullet};<3mm,2.4mm>*{}**@{},
<3.3mm,2.9mm>*{};<1.4mm,1.0mm>*{}**@{-},
<0.9mm,0.5mm>*{};<-0.5mm,-0.8mm>*{}**@{-},
<-0.5mm,-0.8mm>*{};<-0.5mm,-3.5mm>*{}**@{-},
<3mm,2.4mm>*{};<3mm,5.2mm>*{}**@{-},
\end{xy}
 +
\begin{xy}
<3mm,-0.8mm>*{\circ};<3mm,-0.8mm>*{}**@{},
<3mm,-1.3mm>*{};<3mm,-3.5mm>*{}**@{-},
<2.3mm,-0.5mm>*{};<0.3mm,1.2mm>*{}**@{-},
<0.3mm,1.2mm>*{};<2.4mm,2.4mm>*{}**@{-},
<3.5mm,-0.5mm>*{};<7mm,2.4mm>*{}**@{-},
<7mm,2.4mm>*{};<7mm,5.2mm>*{}**@{-},
<3mm,2.4mm>*{\bullet};<3mm,2.4mm>*{}**@{},
<2.7mm,2.9mm>*{};<4.8mm,1.0mm>*{}**@{-},
<5.4mm,0.5mm>*{};<7mm,-0.8mm>*{}**@{-},
<7mm,-0.8mm>*{};<7mm,-3.5mm>*{}**@{-},
<3mm,2.4mm>*{};<3mm,5.2mm>*{}**@{-},
\end{xy}
$$

\vspace{1mm}

                       \end{tabular}

  & \begin{tabular}{l} \vspace{- 3mm}\\
                      smooth formal Poisson\\
                      structure in $\R^n$ \\
                      vanishing at $0$ \vspace{1 mm}\\
                       \vspace{1 mm}

                      \end{tabular}
&
\begin{tabular}{l} \vspace{- 3mm}\\
                    structure, $(\widehat{V\oplus V^*[1]},\eth)$,\\
                      of a smooth dg manifold \\
                       together with an odd\\
                      symplectic form $\om_{odd}$\\
                      on $\widehat{V\oplus V^*[1]}$ such that \\
                      the homological vector\\
                      field $\eth$ is hamiltonian\\
                      and vanishes on $\widehat{0\oplus V^*[1]}$\vspace{1 mm}\\
                       \end{tabular}

 \\
\hline \multicolumn{3}{|c|}{\begin{tabular}{ll} & \\
{\sc
Notations:} &
For a graded vector space $V$, $\hat{V}$ stands for the formal graded manifold\\
& (non-canonically) isomorphic to the formal neighbourhood of $0$ in $V$,\\
& and $T_{\hat{V}}$ stands for
 the tangent bundle on $\hat{V}$.\end{tabular} }\\

\hline

\end{tabular}

\end{center}
\end{table}


\bip

{\bf 1.3. Hertling-Manin's geometry and the $G$-operad.} A
{\em Gerstenhaber}\,
algebra  is, by definition, a graded vector space $V$ together with two
linear maps,
$$
\Ba{rccc}

\circ: & \odot^2 V & \lon & V \\
       & a\ot b    & \lon & a\circ b
\Ea
\ \ \ \ \  , \ \ \ \
\Ba{rccc}

[\, \bullet\, ]: & \odot^2 V & \lon & V[1] \\
       & a\ot b    & \lon & (-1)^{|a|}[a\bullet b]
\Ea
$$
satisfying the identities,
\Bi
\item[(i)] $a\circ (b\circ c) -
(a\circ b)\circ c =0$ (associativity);
\item[(ii)] $[[a\bullet
b]\bullet c]=[a\bullet[b\bullet c]] +
(-1)^{|b||a|+|b|+|a|}[b\bullet[a\bullet c]]$ (Jacobi identity);
\item[(iii)]
 $[(a\circ b)\bullet c]= a\circ [b\bullet c] +    (-1)^{|b|(|c|+1)}[a\bullet c]\circ b$
(Leibniz type identity).
\Ei

The operad whose algebras are Gerstenhaber algebras is often called the $G$-operad. It
has a relatively simple structure, $Free(E)/Ideal$, with $E$ spanned by two corollas,
$$
E=span\left\{ \circ =
 \begin{xy}
 <0mm,0.66mm>*{};<0mm,3mm>*{}**@{-},
 <0.39mm,-0.39mm>*{};<2.2mm,-2.2mm>*{}**@{-},
 <-0.35mm,-0.35mm>*{};<-2.2mm,-2.2mm>*{}**@{-},
 <0mm,0mm>*{\circ};<0mm,0mm>*{}**@{},
 \end{xy}\  , \,
[\, \bullet\, ]=
 \begin{xy}
 <0mm,0.66mm>*{};<0mm,3mm>*{}**@{-},
 <0.39mm,-0.39mm>*{};<2.2mm,-2.2mm>*{}**@{-},
 <-0.35mm,-0.35mm>*{};<-2.2mm,-2.2mm>*{}**@{-},
 <0mm,0mm>*{\bullet};<0mm,0mm>*{}**@{},
 \end{xy}
\right\}
$$
and with engineering rules (i)-(iii). The minimal resolution of the
$G$-operad has been constructed in \cite{GJ} and is often called a
$G_\infty$-operad. The derived category of Gerstenhaber algebras
is equivalent to the category whose objects are  isomorphism
classes of minimal $G_\infty$-structures on graded vector spaces
$V$. Let $(M,*)$ be the formal pointed graded manifold whose
tangent space at the distinguished point is isomorphic to a vector
space $V$, and let us choose an arbitrary torsion-free affine
connection $\nabla$ on $M$. With this choice a structure of
$G_\infty$ algebra on a graded vector space $V$ can be suitable
described as
\Bi
\item a degree 1 smooth vector
field $\eth$ on $M$ satisfying the integrability condition
$[\eth,\eth]=0$ and vanishing at the distinguished point $*$; (if
$\eth$ has zero at $*$ of second order, then the
$G_\infty$-structure is called minimal);
\item  a collection of
homogeneous {\em tensors},
$$
\left\{
\mu_{n_1,\ldots, n_k}: T_M^{\ot n_1} \ot
T_M^{\ot n_2} \ot \ldots \ot   T_M^{\ot n_k} \rar T_M[k+1-n_1-\ldots - n_k]
\right\}_{n_i,k\geq 1,n_i + k\geq 2}
$$
satisfying an infinite tower of quadratic algebraic and differential equations.
The first two floors of this tower read as follows: the data
$\{\mu_{n}\}_{n\geq 1}$
(with $ \mu_1:=Lie_\eth$)
makes the tangent sheaf $T_M$ into a sheaf of $C_{\infty}$
algebras\footnote{$C_{\infty}$ stands
for the minimal resolution of the operad of commutative associative algebras.}
satisfying an ``integrability'' condition,
$$
[\mu_{\bullet},\mu_{\bullet}]_{G_\infty}=Lie_\eth \mu_{\bullet,\bullet}
$$
for a certain bi-differential operator $[\ ,\ ]_{G_\infty}$ whose leading term is just
the usual vector field bracket of values of $\mu_\bullet$.  It is also required that
each tensor $\mu_{\bullet, \ldots, \bullet}:
T_M^{\ot \bullet} \ot \ldots \ot   T_M^{\ot \bullet} \rar T_M$
vanishes if the input contains at least one pure shuffle product,
$$
(v_1\ot \ldots \ot v_k)\bigstar (v_{k+1}\ot\ldots \ot v_n):=
\sum_{{\mathrm Shuffles}\ \sigma\atop
{\mathrm of \ type\ (k,n)}} (-1)^{Koszul(\sigma)}v_{\sigma(1)}\ot \ldots
\ot v_{\sigma(n)},\ \ v_i\in T_M.
$$

\Ei

A change of the connection $\nabla$ alters the tensors
$\mu_{\bullet_1,\ldots, \bullet_k}$, $k\geq 2$, but leaves the
homotopy class of the $G_\infty$-structure on $V$ invariant. \sip

\sip

If the vector space $V$ is concentrated in degree $0$, i.e.
$V\simeq \R^n$, then a $G_\infty$-structure
on $V$ reduces just to a single tensor field $\mu_2: T_M^{\ot 2}\rar T_M$ which makes
 the tangent sheaf
into a sheaf of commutative associative algebras, and satisfies the differential equations,
 $$
[\mu_{2},\mu_{2}]_{G_\infty}=0.
$$
The explicit form for the bracket $[\ , \ ]_{G_\infty}$ can be read off from the $G_\infty$
operad structural equations rather straightforwardly (see \cite{Me} for details),
\Beqrn
  [\mu_2,\mu_2]_{G_\infty}(X,Y,Z,W)&=&
  [\mu_2(X, Y), \mu_2(Z,W)] - \mu_2([\mu_2(X, Y), Z], W)
  - \mu_2(Z,[\mu_2(X,Y), W])
   \\
 && - \mu_2(X, [Y,\mu_2(Z, W)]) -  \mu_2[X,
 \mu_2(Z,W)], Y) \\
 && + \mu_2(X, \mu_2(Z, [Y,W]))  + \mu_2(X, \mu_2([Y,Z], W)) \\
 && + \mu_2([X,Z],\mu_2( Y, W)) + \mu_2([X,W],\mu_2(Y, Z)).
  \Eeqrn
The resulting geometric structure is precisely the one
discovered earlier by Hertling and Manin \cite{HM} in their quest
for a weaker notion of Frobenius manifold; they call it an
$F$-{\em manifold}\, structure on $V$.

\sip

Hertling-Manin's geometric structures arise naturally  in the
theory of singularities \cite{H} and the deformation theory
\cite{Me}.

 \bip

{\bf 1.4. Germs of tensor fields.}

 A ${\mathsf TF}$ {\em bialgebra}  is, by definition, a graded vector space $V$ together
with two linear maps,
$$
\Ba{rccc}
\delta \equiv
 \begin{xy}
 <0mm,-0.55mm>*{};<0mm,-2.5mm>*{}**@{-},
 <0.5mm,0.5mm>*{};<2.2mm,2.2mm>*{}**@{-},
 <-0.48mm,0.48mm>*{};<-2.2mm,2.2mm>*{}**@{-},
 <0mm,0mm>*{\circ};<0mm,0mm>*{}**@{},
 \end{xy}: &  V & \lon & \ot^2 V \\
       & a    & \lon & \sum a_{1}\ot a_{2}
\Ea \ \ \ \ \  , \ \ \ \ \Ba{rccc}
[\, \bullet\, ]\equiv
 \begin{xy}
 <0mm,0.66mm>*{};<0mm,3mm>*{}**@{-},
 <0.39mm,-0.39mm>*{};<2.2mm,-2.2mm>*{}**@{-},
 <-0.35mm,-0.35mm>*{};<-2.2mm,-2.2mm>*{}**@{-},
 <0mm,0mm>*{\bullet};<0mm,0mm>*{}**@{},
 \end{xy}: & \odot^2 V& \lon & V[1] \\
       & a\ot b    & \lon & (-1)^{|a|}[a\bullet b]
\Ea
$$
satisfying the identities,
\Bi
 \item[(i)]
$[[a\bullet b]\bullet c]=[a\bullet[b\bullet c]] +
(-1)^{|b||a|+|b|+|a|}[b\bullet[a\bullet c]]$ (Jacobi identity);
\item[(ii)]
 $\delta [a\bullet b]=\sum a_1\ot [a_2\bullet b] +  [a\bullet b_1]\ot b_2
+ (-1)^{|a||b|+|a|+|b|}(
[b\bullet a_1]\ot a_2  + b_1\ot [b_2\bullet a])$ (Leibniz type
identity).
\Ei

There are obvious versions of the above notion  with $\delta$ taking values
in $\wedge^2 V$ and $\odot^2 V$, i.e.\ with the gene
$
\begin{xy}
 <0mm,-0.55mm>*{};<0mm,-2.5mm>*{}**@{-},
 <0.5mm,0.5mm>*{};<2.2mm,2.2mm>*{}**@{-},
 <-0.48mm,0.48mm>*{};<-2.2mm,2.2mm>*{}**@{-},
 <0mm,0mm>*{\circ};<0mm,0mm>*{}**@{},
 \end{xy}
 $
  realizing either the trivial or sign representations of $\Sigma_2$.

\sip

The dioperad whose algebras are  ${\mathsf TF}$ bialgebras  is denoted by
${\mathsf TF}$. This quadratic dioperad is Koszul so that one can construct its minimal
resolution using the results of \cite{G,GK,Mar}. It turns out that the structure
of ${\mathsf  TF}_\infty$-algebra on a graded vector space $V$ is the same as a pair of
collections of linear maps,
$$
\left\{\mu_{n}: \odot^n V\rar V[1]\right\}_{n\geq 1},
$$
and
$$
\left\{\phi_{n}: \odot^n V\rar V\ot V\right\}_{n\geq 1},
$$
 satisfying a system of
quadratic equations which are best described using a geometric language.
Let $M$ be the formal graded
manifold associated to $V$. If $\{e_\al, \al=1,2,\ldots\}$
is a homogeneous basis of $V$, then the associated dual basis $t^{\al}$,
 $|t^{\al}|=-|e_{\al}|$,
defines a coordinate system on $M$. The collection of tensors
$\{\mu_{n}\}_{n\geq 1}$ can be assembled into a germ, $\eth\in T_M$, of
a degree 1 smooth vector field,
$$
 \eth:= \sum_{n=1}^{\infty}\frac{1}{n!} (-1)^{\epsilon}t^{\al_1}\cdots t^{\al_n}
\mu_{\al_1 \ldots \al_n}^{\ \ \ \ \ \ \be} \frac{\p}{\p t^{\be}}
$$
where
$$
\epsilon=\sum_{k=1}^n|e_{\al_k}|(1+\sum_{i=1}^k|e_{\al_i}|)
$$
the numbers $\mu_{\al_1 \ldots \al_n}^{\ \ \ \ \ \ \be}$ are defined by
$$
\mu_{n}(e_{\al_1}, \ldots , e_{\al_n})=
\sum \mu_{\al_1 \ldots \al_n}^{\ \ \ \ \ \ \be}e_{\be},
$$
and we assume here and throughout the paper  summation over repeated small Greek indices.

\sip

Another collection of linear maps, $\{\phi_n\}$, can be assembled into a smooth germ,
$\phi\in \ot^2 T_M$, of a degree zero contravariant tensor field on $M$,
$$
\phi:= \sum_{n=1}^{\infty}\frac{1}{n!} (-1)^{\epsilon}t^{\al_1}\cdots t^{\al_n}
\phi_{\al_1 \ldots \al_n}^{\ \ \ \ \ \ \be_1\be_2} \frac{\p}{\p t^{\be}}\ot
\frac{\p}{\p t^{\be}}
$$
where
$$
\epsilon=|e_{\be_2}|(|e_{\be_1}| +1)+
\sum_{k=1}^n\sum_{i=1}^k|e_{\al_k}||e_{\al_i}|
$$
and the numbers $\mu_{\al_1 \ldots \al_n}^{\ \ \ \ \ \ \be_1\be_2}$ are defined by
$$
\mu_{n}(e_{\al_1}, \ldots , e_{\al_n})=
\sum \mu_{\al_1 \ldots \al_n}^{\ \ \ \ \ \ \be_1\be_2}e_{\be_1}\ot e_{\be_2}.
$$

\mbox{\bf 1.4.1. Proposition.} {\em The collections of tensors,
$$
\left\{\mu_{n}: \odot^n V\rar V[1]\right\}_{n\geq 1}\ \ \  and \ \ \
\left\{\phi_{n}: \odot^n V\rar V\ot V\right\}_{n\geq 1},
$$
define a structure of $\mathsf TF_\infty$-algebra
on $V$ if and only if the associated smooth vector field $\eth$ and the contravariant
tensor field $\phi$ satisfy the equations,
$$
[\eth,\eth]=0
$$
and
$$
Lie_\eth\phi=0,
$$
where $[\ , \ ]$ stands for the usual bracket of vector fields and $Lie_\eth$
for the Lie derivative along $\eth$.}

\bip

If $V$ is finite dimensional and concentrated in degree zero, then
a $\mathsf TF_\infty$-structure in $V$ is just a germ of a smooth
rank 2 contravariant tensor on $V$ vanishing at $0$.

 \bip

{\bf 1.5. Poisson geometry and the dioperad of Lie
1-bialgebras.}

 A {\em  Lie
1-bialgebra}\,
  is, by definition, a graded vector space $V$ together with two linear maps,
$$
\Ba{rccc}
\delta: &  V & \lon & \wedge^2 V \\
       & a    & \lon & \sum a_{1}\wedge a_{2}
\Ea
\ \ \ \ \  , \ \ \ \
\Ba{rccc}
[\, \bullet\, ]: & \odot^2 V& \lon & V[1] \\
       & a\ot b    & \lon & (-1)^{|a|}[a\bullet b]
\Ea
$$
satisfying the identities,
\Bi
\item[(i)] $(\delta\ot\Id)\delta a
+ \tau (\delta\ot\Id)\delta a+ \tau^2 (\delta\ot\Id)\delta a =0$,
where $\tau$ is the cyclic permutation $(123)$ represented
naturally on $V\ot V \ot V$ (co-Jacobi identity); \item[(ii)]
$[[a\bullet b]\bullet c]=[a\bullet[b\bullet c]] +
(-1)^{|b||a|+|b|+|a|}[b\bullet[a\bullet c]]$ (Jacobi identity);
\item[(iii)]
 $\delta [a\bullet b]=\sum a_1\wedge [a_2\bullet b] - (-1)^{|a_1||a_2|} a_2\wedge
[a_1\bullet b]
+  [a\bullet b_1]\wedge b_2 - (-1)^{|b_1||b_2|}[a\bullet b_2]\wedge b_1$
(Leibniz type identity).
\Ei

The dioperad whose algebras are  Lie 1-bialgebras  is denoted by
$\LieBi$. The superscript $1$ in the notation is used to
emphasize that the two basic operations
$$
\delta =
 \begin{xy}
 <0mm,-0.55mm>*{};<0mm,-2.5mm>*{}**@{-},
 <0.5mm,0.5mm>*{};<2.2mm,2.2mm>*{}**@{-},
 <-0.48mm,0.48mm>*{};<-2.2mm,2.2mm>*{}**@{-},
 <0mm,0mm>*{\circ};<0mm,0mm>*{}**@{},
 \end{xy}
 \ \ \  , \ \ \
[\, \bullet\, ]=
 \begin{xy}
 <0mm,0.66mm>*{};<0mm,3mm>*{}**@{-},
 <0.39mm,-0.39mm>*{};<2.2mm,-2.2mm>*{}**@{-},
 <-0.35mm,-0.35mm>*{};<-2.2mm,-2.2mm>*{}**@{-},
 <0mm,0mm>*{\bullet};<0mm,0mm>*{}**@{},
 \end{xy}
$$
have homogeneities differed  by $1$.

\sip

Similarly one can introduce the notion of {\em Lie $n$-bialgebras}:
coLie algebra structure on $V$ plus Lie algebra structure on
$V[-n]$  plus an obvious Leibniz type identity. Homotopy theory of
Lie $n$-bialgebras splits into two stories, one for $n$ even, and
one for $n$ odd. The even case (more precisely, the case $n=0$)
has been studied by Gan \cite{G}. In this paper we study the odd
case, more precisely, the case $n=1$.

\sip

The dioperad $\LieBi$ is Koszul. Hence one can use the machinery
of \cite{G,GK,Mar} to construct its minimal resolution, the
dioperad $\LieBi_\infty$. The structure of a
$\LieBi_\infty$ algebra on a graded vector space $V$ is a
collection linear maps,
$$
\left\{\mu_{m,n}: \odot^n V\rar \wedge^m V[2-m]\right\}_{m\geq 1, n\geq 1},
$$
 satisfying a system of
quadratic equations which can be described as follows. Let $M$ be the formal graded
manifold associated to $V$. If $\{e_\al, \al=1,2,\ldots\}$
is a homogeneous basis of $V$, then the associated dual basis $t^{\al}$,
 $|t^{\al}|=-|e_{\al}|$,
defines a coordinate system on $M$. For a fixed $m$ the collection of tensors
$\{\mu_{m,n}\}_{n \geq 1}$ can be assembled into a germ, $\Gamma_m\in \wedge^m T_M$, of
a smooth polyvector field (vanishing at $0\in M$),
$$
 \Gamma_m:= \sum_{n=1}^{\infty}\frac{1}{m!n!} (-1)^{\epsilon}t^{\al_1}\cdots t^{\al_n}
\mu_{\al_1 \ldots \al_n}^{\ \ \ \ \ \ \be_1\ldots \be_m} \frac{\p}{\p t^{\be_1}}
\wedge \ldots \wedge  \frac{\p}{\p t^{\be_m}}
$$
where
$$
\epsilon=\sum_{k=1}^n|e_{\al_k}|(2-m+\sum_{i=1}^k|e_{\al_i}|)
+ \sum_{k=1}^n(|e_{\be_k}|+1)\sum_{i=k+1}^n|e_{\be_i}|
$$
and the numbers $\mu_{\al_1 \ldots \al_n}^{\ \ \ \ \ \ \be_1\ldots \be_m}$ are defined by
$$
\mu_{m,n}(e_{\al_1}, \ldots , e_{\al_n})=
\sum \mu_{\al_1 \ldots \al_n}^{\ \ \ \ \ \ \be_1\ldots \be_m}e_{\be_1}\wedge \ldots
\wedge e_{\be_m}.
$$

\bip

{\bf 1.5.1. Proposition.} {\em A collection of tensors, $
\left\{\mu_{m,n}: \odot^n V\rar \wedge^m V[2-m]\right\}_{m\geq 1
, n\geq 1}$, defines a structure of $\LieBi_\infty$-algebra
on $V$ if and only if the associated smooth polyvector field,
$$
\Gamma := \sum_{m\geq 1} \Gamma_m\in \wedge^{\bullet}T_M,
$$
satisfies the equation
$$
[\Gamma, \Gamma]=0,
$$
where $[\ , \ ]$ stands for the Schouten bracket of polyvector fields.}

\bip

In particular, if $V$ is concentrated in degree zero, then
the only non-zero summand in $\Gamma$ is $\Gamma_2\in \wedge^2 T_M$. Hence a
 $\LieBi_\infty$-algebra structure on $\R^n$ is nothing but a germ of a smooth
Poisson structure on  $\R^n$ vanishing at $0$.

\bip

{\bf 1.6. On the  content of the rest.} 
Section 2 is a reminder on PROPs, dioperads and
Koszulness \cite{G,GK,Mar}.
In Sections 3 and 4  we prove Koszulness of the dioperads $\LieBi$ and ${\mathsf TF}$,
apply the machinery reviewed in Section 2 to give  explicit graph descriptions
of their minimal resolutions,   $\LieBi_\infty$ and $\mathsf TF_\infty$,
prove Propositions 1.4.1 and 1.5.1 and introduce and study the notion of
$\LieBi_\infty$ morphisms.
Section 5 is a comment on a geometric description
of algebras over the dioperad of strongly homotopy Lie bialgebras,
 and  their strongly homotopy maps.

\bip

\sip

{\bf 2. PROPs and dioperads \cite{G}.} Let ${\mathsf S_f}$ be the groupoid of finite sets.
It is equivalent to the category whose objects
are natural numbers, $\{m\}_{m\geq 1}$, and morphisms are the
permutation groups $\{\Sigma_m\}_{m\geq 1}$.

\sip

 A PROP
${\cP}$ in the category, ${\mathsf dgVec}$, of differential graded (shortly, dg)
vector spaces
is a functor $\cP: {\mathsf S_f}\times {\mathsf S_f}^{op}\rar {\mathsf dgVec}$
together with natural transformations,
\Beqrn
\circ_{A,B,C}: && \cP(A,B)\ot \cP(B,C) \lon \cP(A,C),\\
\ot_{A,B,C,D}: &&  \cP(A,B)\ot \cP(C,D) \lon \cP(A\ot B,C\ot D)
\Eeqrn
and the distinguished elements $\Id_A\in \cP(A,A)$ and
$s_{A,B}\in \cP(A\ot B, B\ot A)$
satisfying a system of axioms \cite{A} which just mimic the obvious properties of
the following natural transformation,
$$
\cE_V: (m,n) \lon  Hom(V^{\ot n}, V^{\ot m}),
$$
canonically associated with an arbitrary dg space $V$. The latter fundamental example
 is called the {\em endomorphism}\,
PROP of $V$.

\sip

Given a collection of dg $(\Sigma_m, \Sigma_n)$-bimodules,
$E=\{E(m,n)\}_{m,
n\geq 1}$,
one can construct the associated free PROP, $Free(E)$,
 by decorating vertices of all possible
directed graphs with a flow by the elements of $E$ and then taking
the colimit over the graph automorphism group. The composition operation $\circ$
corresponds then to gluing output legs of one graph to the input legs of another graph,
and the tensor product $\ot$ to the disjoint union of graphs. Even for a  small
finite dimensional collection $E$ the resulting free PROP can be a monstrous  infinite
dimensional object. The notion of dioperad was
introduced by Gan \cite{G} as a way to avoid that free PROP ``divergence''.
In the above setup, a free dioperad on $E$ is built
on graphs of genus zero, i.e.\ on trees.

\sip

More precisely, a {\em dioperad}\, $\cP$ consists of data:
\Bi
\item[(i)] a   collection of dg $(\Sigma_m,\Sigma_n)$ bimodules,
$\{\cP(m,n)\}_{m\geq 1,
n\geq 1}$;
\item[(ii)] for each $m_1,n_1,m_2,n_2\geq 1$,
 $i\in \{1,2,\ldots, n_1\}$ and  $j\in \{1,\ldots, n_1\}$ a linear
map
$$
_i\circ_j: \cP(m_1, n_1)\ot \cP(m_2,n_2)\lon \cP(m_1+m_2-1, n_1+n_2-1),
$$
\item[(iii)] a morphism  $e: k\rar \cP(1,1)$ such that the compositions
$$
k\ot \cP(m,n) \stackrel{e\ot Id}{\lon} \cP(1,1)\ot \cP(m,n)
 \stackrel{_1\circ_i}{\lon}\cP(m,n)
$$
and
$$
\cP(m,n)\ot k \stackrel{Id\ot e}{\lon} \cP(m,n)\ot \cP(1,1)
 \stackrel{_j\circ_1}{\lon}\cP(m,n)
$$
are the canonical isomorphisms for all $m,n\geq 1$, $1\leq i\leq m$ and
  $1\leq j\leq n$.
\Ei
These data satisfy associativity and equivariance conditions \cite{G} which can be read off
from
the example of the {\em endomorphism dioperad}\, $\cE nd_V$ with
$\cE nd_V(m,n)=Hom(V^{\ot n},V^{\ot m})$, $e: 1\rar \Id\in Hom(V,V)$, and
the compositions given by
$$
\Ba{cccc}

_i\circ_j: &  \cP(m_1, n_1)\ot \cP(m_2,n_2) & \lon & \cP(m_1+m_2-1, n_1+n_2-1) \\
           &   f\ot g                       & \lon &
(\Id\ot \ldots \ot f\ot \ldots \ot\Id)\sigma (\Id\ot \ldots \ot g\ot \ldots \ot \Id),
\Ea
$$
where $f$ (resp.\ $g$) is at the $j$th (resp.\ $i$th)  place, and
$\sigma$ is the
permutation of the set $I=(1,2, \ldots, n_1+m_2-1)$ swapping the subintervals,
$I_1 \leftrightarrow I_2$ and $I_4\leftrightarrow I_5$, of the unique
order preserving  decomposition, $I=I_1 \sqcup I_2  \sqcup I_3 \sqcup I_4 \sqcup I_5$,
of $I$ into the disjoint union of five intervals of lengths
$|I_1|=i-1$,  $|I_2|=j-1$,  $|I_3|=1$,  $|I_4|=m_2-j$
and  $|I_5|=n_1-i$.

\sip

If $\cP$ is a dioperad, then the collection of   $(\Sigma_m,\Sigma_n)$ bimodules,
$$
\cP^{op}(m,n):=\left(\cP(n,m), {\mathrm transposed\ actions\ of}\
 \Sigma_m\ {\mathrm and }\ \Sigma_n\right),
$$
is naturally a dioperad as well.

\sip

If $\cP$ is a dioperad with $\cP(m,n)$ vanishing for all $m,n$ except
for $(m=1, n\geq 1)$,
then $\cP$ is called an {\em operad}.

\sip

A morphism of dioperads, $F: \cP\rar {\cal Q}$, is a collection of equivariant
linear maps, $F(m,n): \cP(m,n)\rar {\cal Q}(m,n)$, preserving all the structures.
If $\cP$ is a dioperad, then a $\cP$-{\em algebra} is a dg vector space
$V$ together with a morphism, $F:\cP\rar \cE nd_V$, of dioperads.

\sip

We shall consider below only  dioperads $\cP$
with $\cP(m,n)$ being finite
dimensional vector spaces (over a field $k$ of characteristic zero)
for all $m,n$.

\sip

The endomorphism dioperad of the vector space $k[-p]$, $p\in \Z$,  is denoted
by $\langle p \rangle$. Thus $\langle p \rangle(m,n)$ is
$sgn_n^{\ot p}\ot sgn_m^{\ot p} [p(n-m)]$
where $sgn_m$ stands for the one dimensional sign representation of $\Sigma_m$.
Representations of the dioperad $\cP\langle p \rangle:=\cP\ot \langle p \rangle$
in a vector space $V$
are the same as representations of the dioperad $\cP$ in $V[p]$.

\sip

If $\cP$ is a dioperad, then $\Lambda\cP:=\{sgn_m\ot\cP(m,n)[2-m-n]\ot sgn_n\}$
 and
$\Lambda^{-1}\cP:= =\{sgn_m\ot\cP(m,n)[m+n-2]\ot sgn_n\}$ are also dioperads.

\bip

{\bf 2.1. Cobar dual.} If $T$ is a directed (i.e.\ provided with a flow
which we always assume in our pictures to go from the bottom to the top)
tree, we denote by
\Bi
\item $Vert(T)$ the set of all vertices,
\item $edge(T)$ the set
of internal edges; $\det(T):= \wedge^{|edge(T)|}{\mathrm span}_k(edge(T))$;
\item   $Edge(T)$ is  the set of all edges, i.e.\
$$
Edge(T):=edge(T)\sqcup\{\mathrm input\ legs\ (leaves)\}\sqcup
\{\mathrm output\ legs\ (roots)\};
$$
 ${\mathrm Det}(T):= \wedge^{|Edge(T)|}{\mathrm span}_k(Edge(T))$;
\item $Out(v)$ (resp.\ $ In(v)$) the set of outgoing (resp.\, incoming)
edges at a vertex $v\in Vert(V)$.
\Ei
An $(m,n)$-tree is a tree $T$ with $n$ input legs labeled by
the set $[n]=\{1,\ldots,n\}$
and $m$ output legs  labeled by the set $[m]=\{1,\ldots,m\}$.
A tree $T$ is called {\em trivalent}\, if $|Out(v)\sqcup In(v)|=3$
for all $v\in Vert(T)$.

\sip

Let $E=\{E(m,n)\}_{m,n\geq 1}$ be a collection of finite
dimensional $(\Sigma_m,\Sigma_n)$ bimodules with $E_{1,1}=0$. For a pair of finite
sets, $I,J\in Objects(\mathsf S_f)$, with $|I|=m$ and $|J|=n$, one defines
$$
E(I,J):= {Hom}_{\mathsf S_f}([m],I)\times_{\Sigma_m} E(m,n))\times_{\Sigma_n}
 Hom_{\sf S_f}(J,[n]).
$$

\sip

 The {\em free
dioperad}, $Free(E)$, generated by $E$ is defined by
$$
Free(E)(m,n):= \bigoplus_{(m,n)-{\mathrm trees}\ T}\ E(T),
$$
where
$$
E(T):=\bigotimes_{v\in Vert(T)}E(Out(v), In(v)),
$$
and the compositions $_i\circ_j$ are given by grafting the $j$th root of one tree into
$i$th leaf of another tree, and then taking  the ``unordered'' tensor product
\cite{MSS}  over the
set of vertices of the resulting tree.

\sip

Let $\cP=\{\cP(m,n)\}_{m,n\geq 1}$ be a collection of graded
$(\Sigma_m,\Sigma_n)$ bimodules.
We denote by $\bar{\cP}$ the collection
$\{\bar{\cP}(m,n)\}_{m,n\geq 1}$ given by $\bar{\cP}(m,n):= {\cP}(m,n)$
for $m+n\geq 3$
and $\bar{\cP}(1,1)=0$. The collection of dual vector spaces,
$\bar{\cP}^*=\{\bar{\cP}(m,n)^*\}_{m,n\geq 1}$, is naturally a collection of
$(\Sigma_m,\Sigma_n)$-bimodules with the transposed actions. We also set
  $\cP^{\vee}=\{\bar{\cP}(m,n)^\vee:= sgn_m\ot\bar{\cP}(m,n)^*\ot sgn_n\}$.

\sip

Let $\cP$ be a graded dioperad with zero differential.
The {\em cobar dual}\,
of $\cP$ is the dg dioperad ${\mathbf D}\cP$ defined by

(i) as a dioperad of graded vector spaces,
${\mathbf D}\cP= \Lambda^{-1}Free(\bar{\cP}^*[-1])=
Free(\Lambda^{-1}\bar{\cP}^*[-1])$;

(ii) as a complex,  ${\mathbf D}\cP$ is non-positively graded,  ${\mathbf D}\cP(m,n) =
\sum_{i=0}^{m+n-3}{\mathbf D}\cP^{-i}(m,n)$ with the differential
 given by dualizations of the compositions $_\bullet\circ_\bullet$
and edge contractions \cite{G,GK},
$$
\Ba{ccccccccc}
{\mathbf D}\cP^{3-m-n}(m,n) &\hspace{-2mm}\stackrel{d}{\rar}\hspace{-2mm}&
{\mathbf D}\cP^{4-m-n}(m,n)  &\hspace{-2mm}\stackrel{d}{\rar}\hspace{-2mm}&
{\mathbf D}\cP^{3-m-n}(m,n)   &\hspace{-2mm}\stackrel{d}{\rar}\hspace{-2mm}&
\ldots
 &\hspace{-2mm}\stackrel{d}{\rar}\hspace{-2mm}&
{\mathbf D}\cP^{0}(m,n)\\
|| & & || && || &&&& || \\
\bar{\cP}^\vee(m,n)
 &\hspace{-2mm}\stackrel{d}{\rar}\hspace{-2mm}&
\underset{|edge(T)|=1}{\bigoplus}\bar{\cP}^*\ot {\mathrm Det}(T)
 &\hspace{-2mm}\stackrel{d}{\rar}\hspace{-2mm}&
\underset{|edge(T)|=2}{\bigoplus}\bar{\cP}^*\ot {\mathrm Det}(T)
 &\hspace{-2mm}\stackrel{d}{\rar}\hspace{-2mm}&
\ldots
 &\hspace{-2mm}\stackrel{d}{\rar}\hspace{-2mm}&
\underset{{\mathrm trivalent}\atop {\mathrm trees}\ T}
{\bigoplus}\bar{\cP}^*\ot {\mathrm Det}(T)
\Ea
$$
where the sums are taken over $(m,n)$-trees.

\bip

{\bf 2.2. Remark.} The vector space ${\mathbf D}\cP$ is bigraded: one grading comes
from the grading of $\cP$ as a vector space and another one from trees as in (ii) just above.
The differential preserves the first grading and increases by 1 the second one.
The $\Z$-grading of ${\mathbf D}\cP$ is always understood to be the associated
total grading. In particular, $deg_{{\mathbf D}\cP}\bar{\cP}^\vee(m,n)=
deg_{\mathsf Vect}(\bar{\cP}^\vee(m,n)[m+n-3])$.

\bip
{\bf 2.3. Koszul dioperads.} A {\em quadratic dioperad}\, is a dioperad $\cP$
of the form
$$
\cP= \frac{Free(E)}{Ideal <R>},
$$
where $E=\{E(m,n)\}$ is a collection of finite dimensional
$(\Sigma_m,\Sigma_n)$-bimodules
with $E(m,n)=0$ for $(m,n)\neq (1,2), (2,1)$, and the $Ideal$ in
$Free(E)$ is generated by a collection, $R$, of three
sub-bimodules $R(1,2)\subset Free(E)(1,2)$,
 $R(2,1)\subset Free(E)(2,1)$ and  $R(2,2)\subset Free(E)(2,2)$.
 The {\em quadratic dual}\, dioperad, $\cP^!$, is then defined by
$$
\cP^!= \frac{Free(E^\vee)}{Ideal <R^{\perp}>},
$$
where  $R^{\perp}$ is the collection of the three sub-bimodules
 $R^\perp(i,j)\subset Free(E^\vee)(i,j)$ which are annihilators
of  $R(i,j)$, $(i,j)=(1,2), (2,2), (2,1)$.

\sip

Clearly, ${\mathbf D}\cP^0= Free(E^\vee)$ so that there is a natural
epimorphism
$$
{\mathbf D}\cP^0 \lon \cP^!.
$$
Its kernel is precisely $\Img d( {\mathbf D}\cP^{-1})$.
Hence $H^0({\mathbf D}\cP)=  \cP^!$. The quadratic operad $\cP$
is called {\em Koszul}\, if the above morphism is a quasi-isomorphism, i.e.\
$H^{i}({\mathbf D}\cP)=0$ for all $i< 0$. In that case the operad
${\mathbf D}\cP^!$ provides us with a {\em minimal resolution}\, of the operad
$\cP$ and is often denoted by $\cP_\infty$. Algebras over $\cP_\infty$ are
often called {\em strong homotopy}\, $\cP$-algebras; their most important
property is that they can be transferred via quasi-isomorphisms of complexes
\cite{Mar2}.

\bip

{\bf 2.4. Koszulness criterion.} An $(m,n)$-tree $T$ is called {\em reduced}\,
if each vertex has
\Bi
\item either an outgoing root or at least two outgoing internal edges, and/or
\item either an incoming leaf or at least two incoming internal edges.
\Ei
For a collection, $E=\{E_{m,n}\}_{m,n\geq 1}$, of $(\Sigma_m,\Sigma_n)$-bimodules
define another collection of  $(\Sigma_m,\Sigma_n)$-bimodules as follows,
$$
\underline{Free}(E)(m,n):= \bigoplus_{{\mathrm reduced}\atop
(m,n)-{\mathrm trees}} E(T).
$$

\sip

Let $\cP$ be a quadratic dioperad, i.e.\
\mbox{$\cP=Free(E)/Ideal <\hspace{-1mm}R\hspace{-1mm}>$}

for some
generators $E=\{E(1,2), E(2,1)\}$ and relations $R=\{R(1,3), R(2,2), R(3,1)\}$.
With $\cP$ one can canonically associate two quadratic operads, $\cP_L$ and $\cP_R$,
such that
$$
\cP_L= \frac{Free(E(1,2))}{Ideal < R(1,3)>}, \ \ \
\cP_R^{op} =  \frac{Free(E(2,1))}{Ideal <R(2,1)>}.
$$
Let us denote by $\cP_L\diamond \cP_R^{op}$ the collection of  $(\Sigma_m,\Sigma_n)$-bimodules
given by
$$
\cP_L\diamond \cP_R^{op}(m,n):=\left\{
\Ba{cl}
\cP_L(1,n) & {\mathrm if}\ m=1, n\geq 1; \\
\cP_L^{op}(m,1) & {\mathrm if}\ n=1, m\geq 1;\\
0             &{\mathrm otherwise}.
\Ea
\right.
$$

\bip

{\bf 2.4.1. Theorem \cite{G,Mar, MV}}. {\em A quadratic dioperad $\cP$ is Koszul
if the operads $\cP_L$ and $\cP_R$ are Koszul and
$$
\cP(i,j)=\underline{Free}(\cP_L\diamond \cP_R^{op})(i,j)
$$
for $(i,j)=(1,3), (2,2), (3,1)$. Moreover, in this case
$\cP(m,n)=\underline{Free}(\cP_L\diamond \cP_R^{op})(m,n)$ for all\,  $m,n\geq 1$. }

\bip

{\bf 3. A minimal resolution of  $\LieBi$}.
First we present a graph description of the dioperad  $\LieBi$; it will pay off when
we discuss  $\LieBi_\infty$. By definition (see Sect.\ 1.5),  $\LieBi$ is a
quadratic dioperad,
$$
 \LieBi= \frac{ Free(E)}{Ideal <R>}\ ,
$$
where
\Bi
\item[(i)] $E(2,1):= sgn_2\ot 1_1$ and  $E(1,2):= 1_1\ot 1_2[-1]$, where $1_n$
stands for the one dimensional trivial representation of $\Sigma_n$; let
$\delta\in E(2,1)$ and $[\, \bullet \, ]\in E(1,2)$ be basis vectors; we can represent both
as directed\footnote{In all our graphs the direction of edges is chosen to go from the
bottom to the top.}
plane  corollas,
$$
\delta =
 \begin{xy}
 <0mm,-0.55mm>*{};<0mm,-2.5mm>*{}**@{-},
 <0.5mm,0.5mm>*{};<2.2mm,2.2mm>*{}**@{-},
 <-0.48mm,0.48mm>*{};<-2.2mm,2.2mm>*{}**@{-},
 <0mm,0mm>*{\circ};<0mm,0mm>*{}**@{},
 <0mm,-0.55mm>*{};<0mm,-3.8mm>*{_1}**@{},
 <0.5mm,0.5mm>*{};<2.7mm,2.8mm>*{^2}**@{},
 <-0.48mm,0.48mm>*{};<-2.7mm,2.8mm>*{^1}**@{},
 \end{xy}
 \ \ \ \ \  , \ \ \ \ \
[\, \bullet \, ] =
 \begin{xy}
 <0mm,0.66mm>*{};<0mm,3mm>*{}**@{-},
 <0.39mm,-0.39mm>*{};<2.2mm,-2.2mm>*{}**@{-},
 <-0.35mm,-0.35mm>*{};<-2.2mm,-2.2mm>*{}**@{-},
 <0mm,0mm>*{\bullet};<0mm,0mm>*{}**@{},
   <0mm,0.66mm>*{};<0mm,3.4mm>*{^1}**@{},
   <0.39mm,-0.39mm>*{};<2.9mm,-4mm>*{^2}**@{},
   <-0.35mm,-0.35mm>*{};<-2.8mm,-4mm>*{^1}**@{},
\end{xy}
$$
with the following symmetries,
$$
\begin{xy}
 <0mm,-0.55mm>*{};<0mm,-2.5mm>*{}**@{-},
 <0.5mm,0.5mm>*{};<2.2mm,2.2mm>*{}**@{-},
 <-0.48mm,0.48mm>*{};<-2.2mm,2.2mm>*{}**@{-},
 <0mm,0mm>*{\circ};<0mm,0mm>*{}**@{},
   <0mm,-0.55mm>*{};<0mm,-3.8mm>*{_1}**@{},
   <0.5mm,0.5mm>*{};<2.7mm,2.8mm>*{^2}**@{},
   <-0.48mm,0.48mm>*{};<-2.7mm,2.8mm>*{^1}**@{},
 \end{xy}
 = -\,
 \begin{xy}
 <0mm,-0.55mm>*{};<0mm,-2.5mm>*{}**@{-},
 <0.5mm,0.5mm>*{};<2.2mm,2.2mm>*{}**@{-},
 <-0.48mm,0.48mm>*{};<-2.2mm,2.2mm>*{}**@{-},
 <0mm,0mm>*{\circ};<0mm,0mm>*{}**@{},
   <0mm,-0.55mm>*{};<0mm,-3.8mm>*{_1}**@{},
   <0.5mm,0.5mm>*{};<2.7mm,2.8mm>*{^1}**@{},
   <-0.48mm,0.48mm>*{};<-2.7mm,2.8mm>*{^2}**@{},
 \end{xy}
 \ \ \ \ \  , \ \ \ \ \
\begin{xy}
 <0mm,0.66mm>*{};<0mm,3mm>*{}**@{-},
 <0.39mm,-0.39mm>*{};<2.2mm,-2.2mm>*{}**@{-},
 <-0.35mm,-0.35mm>*{};<-2.2mm,-2.2mm>*{}**@{-},
 <0mm,0mm>*{\bullet};<0mm,0mm>*{}**@{},
   <0mm,0.66mm>*{};<0mm,3.4mm>*{^1}**@{},
   <0.39mm,-0.39mm>*{};<2.9mm,-4mm>*{^2}**@{},
   <-0.35mm,-0.35mm>*{};<-2.8mm,-4mm>*{^1}**@{},
\end{xy}
 =
 \begin{xy}
 <0mm,0.66mm>*{};<0mm,3mm>*{}**@{-},
 <0.39mm,-0.39mm>*{};<2.2mm,-2.2mm>*{}**@{-},
 <-0.35mm,-0.35mm>*{};<-2.2mm,-2.2mm>*{}**@{-},
 <0mm,0mm>*{\bullet};<0mm,0mm>*{}**@{},
   <0mm,0.66mm>*{};<0mm,3.4mm>*{^1}**@{},
   <0.39mm,-0.39mm>*{};<2.9mm,-4mm>*{^1}**@{},
   <-0.35mm,-0.35mm>*{};<-2.8mm,-4mm>*{^2}**@{},
\end{xy}\ \  ;
$$
\item[(ii)] the relations $R$ are generated by the following elements,
\Beqrn
\begin{xy}
 <0mm,0mm>*{\circ};<0mm,0mm>*{}**@{},
 <0mm,-0.49mm>*{};<0mm,-3.0mm>*{}**@{-},
 <0.49mm,0.49mm>*{};<1.9mm,1.9mm>*{}**@{-},
 <-0.5mm,0.5mm>*{};<-1.9mm,1.9mm>*{}**@{-},
 <-2.3mm,2.3mm>*{\circ};<-2.3mm,2.3mm>*{}**@{},
 <-1.8mm,2.8mm>*{};<0mm,4.9mm>*{}**@{-},
 <-2.8mm,2.9mm>*{};<-4.6mm,4.9mm>*{}**@{-},
   <0.49mm,0.49mm>*{};<2.7mm,2.3mm>*{^3}**@{},
   <-1.8mm,2.8mm>*{};<0.4mm,5.3mm>*{^2}**@{},
   <-2.8mm,2.9mm>*{};<-5.1mm,5.3mm>*{^1}**@{},
 \end{xy}
\ + \
\begin{xy}
 <0mm,0mm>*{\circ};<0mm,0mm>*{}**@{},
 <0mm,-0.49mm>*{};<0mm,-3.0mm>*{}**@{-},
 <0.49mm,0.49mm>*{};<1.9mm,1.9mm>*{}**@{-},
 <-0.5mm,0.5mm>*{};<-1.9mm,1.9mm>*{}**@{-},
 <-2.3mm,2.3mm>*{\circ};<-2.3mm,2.3mm>*{}**@{},
 <-1.8mm,2.8mm>*{};<0mm,4.9mm>*{}**@{-},
 <-2.8mm,2.9mm>*{};<-4.6mm,4.9mm>*{}**@{-},
   <0.49mm,0.49mm>*{};<2.7mm,2.3mm>*{^2}**@{},
   <-1.8mm,2.8mm>*{};<0.4mm,5.3mm>*{^1}**@{},
   <-2.8mm,2.9mm>*{};<-5.1mm,5.3mm>*{^3}**@{},
 \end{xy}
\ + \
\begin{xy}
 <0mm,0mm>*{\circ};<0mm,0mm>*{}**@{},
 <0mm,-0.49mm>*{};<0mm,-3.0mm>*{}**@{-},
 <0.49mm,0.49mm>*{};<1.9mm,1.9mm>*{}**@{-},
 <-0.5mm,0.5mm>*{};<-1.9mm,1.9mm>*{}**@{-},
 <-2.3mm,2.3mm>*{\circ};<-2.3mm,2.3mm>*{}**@{},
 <-1.8mm,2.8mm>*{};<0mm,4.9mm>*{}**@{-},
 <-2.8mm,2.9mm>*{};<-4.6mm,4.9mm>*{}**@{-},
   <0.49mm,0.49mm>*{};<2.7mm,2.3mm>*{^1}**@{},
   <-1.8mm,2.8mm>*{};<0.4mm,5.3mm>*{^3}**@{},
   <-2.8mm,2.9mm>*{};<-5.1mm,5.3mm>*{^2}**@{},
 \end{xy}
\ \ \ &\in& Free(E)(3,1)\\
&&\\
 \begin{xy}
 <0mm,0mm>*{\bullet};<0mm,0mm>*{}**@{},
 <0mm,0.69mm>*{};<0mm,3.0mm>*{}**@{-},
 <0.39mm,-0.39mm>*{};<2.4mm,-2.4mm>*{}**@{-},
 <-0.35mm,-0.35mm>*{};<-1.9mm,-1.9mm>*{}**@{-},
 <-2.4mm,-2.4mm>*{\bullet};<-2.4mm,-2.4mm>*{}**@{},
 <-2.0mm,-2.8mm>*{};<0mm,-4.9mm>*{}**@{-},
 <-2.8mm,-2.9mm>*{};<-4.7mm,-4.9mm>*{}**@{-},
    <0.39mm,-0.39mm>*{};<3.3mm,-4.0mm>*{^3}**@{},
    <-2.0mm,-2.8mm>*{};<0.5mm,-6.7mm>*{^2}**@{},
    <-2.8mm,-2.9mm>*{};<-5.2mm,-6.7mm>*{^1}**@{},
 \end{xy}
\ + \
 \begin{xy}
 <0mm,0mm>*{\bullet};<0mm,0mm>*{}**@{},
 <0mm,0.69mm>*{};<0mm,3.0mm>*{}**@{-},
 <0.39mm,-0.39mm>*{};<2.4mm,-2.4mm>*{}**@{-},
 <-0.35mm,-0.35mm>*{};<-1.9mm,-1.9mm>*{}**@{-},
 <-2.4mm,-2.4mm>*{\bullet};<-2.4mm,-2.4mm>*{}**@{},
 <-2.0mm,-2.8mm>*{};<0mm,-4.9mm>*{}**@{-},
 <-2.8mm,-2.9mm>*{};<-4.7mm,-4.9mm>*{}**@{-},
    <0.39mm,-0.39mm>*{};<3.3mm,-4.0mm>*{^2}**@{},
    <-2.0mm,-2.8mm>*{};<0.5mm,-6.7mm>*{^1}**@{},
    <-2.8mm,-2.9mm>*{};<-5.2mm,-6.7mm>*{^3}**@{},
 \end{xy}
\ + \
 \begin{xy}
 <0mm,0mm>*{\bullet};<0mm,0mm>*{}**@{},
 <0mm,0.69mm>*{};<0mm,3.0mm>*{}**@{-},
 <0.39mm,-0.39mm>*{};<2.4mm,-2.4mm>*{}**@{-},
 <-0.35mm,-0.35mm>*{};<-1.9mm,-1.9mm>*{}**@{-},
 <-2.4mm,-2.4mm>*{\bullet};<-2.4mm,-2.4mm>*{}**@{},
 <-2.0mm,-2.8mm>*{};<0mm,-4.9mm>*{}**@{-},
 <-2.8mm,-2.9mm>*{};<-4.7mm,-4.9mm>*{}**@{-},
    <0.39mm,-0.39mm>*{};<3.3mm,-4.0mm>*{^1}**@{},
    <-2.0mm,-2.8mm>*{};<0.5mm,-6.7mm>*{^3}**@{},
    <-2.8mm,-2.9mm>*{};<-5.2mm,-6.7mm>*{^2}**@{},
 \end{xy}
\  \ \ &\in & Free(E)(1,3)  \\
&&\\
 \begin{xy}
 <0mm,2.47mm>*{};<0mm,-0.5mm>*{}**@{-},
 <0.5mm,3.5mm>*{};<2.2mm,5.2mm>*{}**@{-},
 <-0.48mm,3.48mm>*{};<-2.2mm,5.2mm>*{}**@{-},
 <0mm,3mm>*{\circ};<0mm,3mm>*{}**@{},
  <0mm,-0.8mm>*{\bullet};<0mm,-0.8mm>*{}**@{},
<0mm,-0.8mm>*{};<-2.2mm,-3.5mm>*{}**@{-},
 <0mm,-0.8mm>*{};<2.2mm,-3.5mm>*{}**@{-},
     <0.5mm,3.5mm>*{};<2.8mm,5.7mm>*{^2}**@{},
     <-0.48mm,3.48mm>*{};<-2.8mm,5.7mm>*{^1}**@{},
   <0mm,-0.8mm>*{};<-2.7mm,-5.2mm>*{^1}**@{},
   <0mm,-0.8mm>*{};<2.7mm,-5.2mm>*{^2}**@{},
\end{xy}
\  - \
\begin{xy}
 <0mm,-1.3mm>*{};<0mm,-3.5mm>*{}**@{-},
 <0.38mm,-0.2mm>*{};<2.2mm,2.2mm>*{}**@{-},
 <-0.38mm,-0.2mm>*{};<-2.2mm,2.2mm>*{}**@{-},
<0mm,-0.8mm>*{\circ};<0mm,0.8mm>*{}**@{},
 <2.4mm,2.4mm>*{\bullet};<2.4mm,2.4mm>*{}**@{},
 <2.5mm,2.3mm>*{};<4.4mm,-0.8mm>*{}**@{-},
 <2.4mm,2.5mm>*{};<2.4mm,5.2mm>*{}**@{-},
     <0mm,-1.3mm>*{};<0mm,-5.3mm>*{^1}**@{},
     <2.5mm,2.3mm>*{};<5.1mm,-2.6mm>*{^2}**@{},
    <2.4mm,2.5mm>*{};<2.4mm,5.7mm>*{^2}**@{},
    <-0.38mm,-0.2mm>*{};<-2.8mm,2.5mm>*{^1}**@{},
    \end{xy}
\  + \
\begin{xy}
 <0mm,-1.3mm>*{};<0mm,-3.5mm>*{}**@{-},
 <0.38mm,-0.2mm>*{};<2.2mm,2.2mm>*{}**@{-},
 <-0.38mm,-0.2mm>*{};<-2.2mm,2.2mm>*{}**@{-},
<0mm,-0.8mm>*{\circ};<0mm,0.8mm>*{}**@{},
 <2.4mm,2.4mm>*{\bullet};<2.4mm,2.4mm>*{}**@{},
 <2.5mm,2.3mm>*{};<4.4mm,-0.8mm>*{}**@{-},
 <2.4mm,2.5mm>*{};<2.4mm,5.2mm>*{}**@{-},
     <0mm,-1.3mm>*{};<0mm,-5.3mm>*{^1}**@{},
     <2.5mm,2.3mm>*{};<5.1mm,-2.6mm>*{^2}**@{},
    <2.4mm,2.5mm>*{};<2.4mm,5.7mm>*{^1}**@{},
    <-0.38mm,-0.2mm>*{};<-2.8mm,2.5mm>*{^2}**@{},
    \end{xy}
\  - \
\begin{xy}
 <0mm,-1.3mm>*{};<0mm,-3.5mm>*{}**@{-},
 <0.38mm,-0.2mm>*{};<2.2mm,2.2mm>*{}**@{-},
 <-0.38mm,-0.2mm>*{};<-2.2mm,2.2mm>*{}**@{-},
<0mm,-0.8mm>*{\circ};<0mm,0.8mm>*{}**@{},
 <2.4mm,2.4mm>*{\bullet};<2.4mm,2.4mm>*{}**@{},
 <2.5mm,2.3mm>*{};<4.4mm,-0.8mm>*{}**@{-},
 <2.4mm,2.5mm>*{};<2.4mm,5.2mm>*{}**@{-},
     <0mm,-1.3mm>*{};<0mm,-5.3mm>*{^2}**@{},
     <2.5mm,2.3mm>*{};<5.1mm,-2.6mm>*{^1}**@{},
    <2.4mm,2.5mm>*{};<2.4mm,5.7mm>*{^2}**@{},
    <-0.38mm,-0.2mm>*{};<-2.8mm,2.5mm>*{^1}**@{},
    \end{xy}
\  + \
\begin{xy}
 <0mm,-1.3mm>*{};<0mm,-3.5mm>*{}**@{-},
 <0.38mm,-0.2mm>*{};<2.2mm,2.2mm>*{}**@{-},
 <-0.38mm,-0.2mm>*{};<-2.2mm,2.2mm>*{}**@{-},
<0mm,-0.8mm>*{\circ};<0mm,0.8mm>*{}**@{},
 <2.4mm,2.4mm>*{\bullet};<2.4mm,2.4mm>*{}**@{},
 <2.5mm,2.3mm>*{};<4.4mm,-0.8mm>*{}**@{-},
 <2.4mm,2.5mm>*{};<2.4mm,5.2mm>*{}**@{-},
     <0mm,-1.3mm>*{};<0mm,-5.3mm>*{^2}**@{},
     <2.5mm,2.3mm>*{};<5.1mm,-2.6mm>*{^1}**@{},
    <2.4mm,2.5mm>*{};<2.4mm,5.7mm>*{^1}**@{},
    <-0.38mm,-0.2mm>*{};<-2.8mm,2.5mm>*{^2}**@{},
    \end{xy}
\ \ \ &\in& Free(E)(2,2).
\Eeqrn
\Ei
\bip

{\bf 3.1. Proposition.} {\em  $\LieBi$ is Koszul.}

\sip

\Proof We have $\LieBi_L={\mathsf Lie}\ot \{\mathsf 1\}$ and $\LieBi_R={\mathsf Lie}$,
where $\mathsf Lie$ stands for the
operad of Lie algebras and
$$
\{\mathsf m\}:=\left\{\{m\}(n):= sgn_n^{\ot m}[m(n-1)]\right\}_{n\geq 1}
$$
for the endomorphism operad of $k[-m]$.
As ${\mathsf Lie}$ is Koszul \cite{GK}, both the operads $\LieBi_L$ and $\LieBi
_R$ are Koszul as well.
Next, a straightforward analysis of all calculational schemes in
$\LieBi$ represented by
directed
trivalent $(i,j)$-trees with $i+j=5$ shows that they generate no new relations so that
$$
\LieBi(i,j)=\underline{Free}\left({\mathsf Lie}\ot\{\mathsf 1\}
\diamond {\mathsf Lie}^{op}\right)(i,j)
$$
for $(i,j)=(1,3), (2,2), (3,1)$. Hence by Theorem~2.4.1, the
dioperad $\LieBi$ is Koszul.
\hfill $\Box$

\bip

Proposition 1.5.1 is a straightforward corollary of the following

\bip

{\bf 3.2. Theorem.}{\em  The minimal resolution,
$\LieBi_\infty$, of the dioperad
 $\LieBi$ can be described as follows.

(i) As a dioperad of graded vector spaces, $\LieBi_\infty=Free(E)$,
where the collection, $E=\{E(m,n)\}$, of one dimensional
$(\Sigma_m,\Sigma_n)$-modules is given by
$$
E(m,n):=\left\{
\Ba{cl}
sgn_m\ot 1_n[m-2] & {\mathrm if}\ m+n\geq 3; \\
0             &{\mathrm otherwise}.
\Ea
\right.
$$

(ii) If we represent a basis element of $E(m,n)$ by the unique (up to a sign)
planar $(m,n)$-corolla,
$$
 
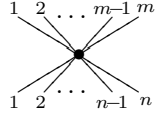
\begin{xy}
 <0mm,0mm>*{\bullet};<0mm,0mm>*{}**@{},
 <0mm,0mm>*{};<-8mm,5mm>*{}**@{-},
 <0mm,0mm>*{};<-4.5mm,5mm>*{}**@{-},
 <0mm,0mm>*{};<-1mm,5mm>*{\ldots}**@{},
 <0mm,0mm>*{};<4.5mm,5mm>*{}**@{-},
 <0mm,0mm>*{};<8mm,5mm>*{}**@{-},
   <0mm,0mm>*{};<-8.5mm,5.5mm>*{^1}**@{},
   <0mm,0mm>*{};<-5mm,5.5mm>*{^2}**@{},
   <0mm,0mm>*{};<4.5mm,5.5mm>*{^{m\hspace{-0.5mm}-\hspace{-0.5mm}1}}**@{},
   <0mm,0mm>*{};<9.0mm,5.5mm>*{^m}**@{},
 <0mm,0mm>*{};<-8mm,-5mm>*{}**@{-},
 <0mm,0mm>*{};<-4.5mm,-5mm>*{}**@{-},
 <0mm,0mm>*{};<-1mm,-5mm>*{\ldots}**@{},
 <0mm,0mm>*{};<4.5mm,-5mm>*{}**@{-},
 <0mm,0mm>*{};<8mm,-5mm>*{}**@{-},
   <0mm,0mm>*{};<-8.5mm,-6.9mm>*{^1}**@{},
   <0mm,0mm>*{};<-5mm,-6.9mm>*{^2}**@{},
   <0mm,0mm>*{};<4.5mm,-6.9mm>*{^{n\hspace{-0.5mm}-\hspace{-0.5mm}1}}**@{},
   <0mm,0mm>*{};<9.0mm,-6.9mm>*{^n}**@{},
 \end{xy}
$$
with skew-symmetric outgoing legs and symmetric ingoing legs, then the differential
$d$ is given on generators by
$$
d \begin{xy}
 <0mm,0mm>*{\bullet};<0mm,0mm>*{}**@{},
 <0mm,0mm>*{};<-8mm,5mm>*{}**@{-},
 <0mm,0mm>*{};<-4.5mm,5mm>*{}**@{-},
 <0mm,0mm>*{};<-1mm,5mm>*{\ldots}**@{},
 <0mm,0mm>*{};<4.5mm,5mm>*{}**@{-},
 <0mm,0mm>*{};<8mm,5mm>*{}**@{-},
   <0mm,0mm>*{};<-8.5mm,5.5mm>*{^1}**@{},
   <0mm,0mm>*{};<-5mm,5.5mm>*{^2}**@{},
   <0mm,0mm>*{};<4.5mm,5.5mm>*{^{m\hspace{-0.5mm}-\hspace{-0.5mm}1}}**@{},
   <0mm,0mm>*{};<9.0mm,5.5mm>*{^m}**@{},
 <0mm,0mm>*{};<-8mm,-5mm>*{}**@{-},
 <0mm,0mm>*{};<-4.5mm,-5mm>*{}**@{-},
 <0mm,0mm>*{};<-1mm,-5mm>*{\ldots}**@{},
 <0mm,0mm>*{};<4.5mm,-5mm>*{}**@{-},
 <0mm,0mm>*{};<8mm,-5mm>*{}**@{-},
   <0mm,0mm>*{};<-8.5mm,-6.9mm>*{^1}**@{},
   <0mm,0mm>*{};<-5mm,-6.9mm>*{^2}**@{},
   <0mm,0mm>*{};<4.5mm,-6.9mm>*{^{n\hspace{-0.5mm}-\hspace{-0.5mm}1}}**@{},
   <0mm,0mm>*{};<9.0mm,-6.9mm>*{^n}**@{},
 \end{xy}
\ \ = \ \
 \sum_{I_1\sqcup I_2=(1,\ldots,m)\atop {J_1\sqcup J_2=(1,\ldots,n)\atop
 {|I_1|\geq 0, |I_2|\geq 1 \atop
 |J_1|\geq 1, |J_2|\geq 0}}
}\hspace{0mm}
(-1)^{\sigma(I_1\sqcup I_2) + |I_1||I_2|}
 \begin{xy}
 <0mm,0mm>*{\bullet};<0mm,0mm>*{}**@{},
 <0mm,0mm>*{};<-8mm,5mm>*{}**@{-},
 <0mm,0mm>*{};<-4.5mm,5mm>*{}**@{-},
 <0mm,0mm>*{};<0mm,5mm>*{\ldots}**@{},
 <0mm,0mm>*{};<4.5mm,5mm>*{}**@{-},
 <0mm,0mm>*{};<13mm,5mm>*{}**@{-},
     <0mm,0mm>*{};<-2mm,7mm>*{\overbrace{\ \ \ \ \ \ \ \ \ \ \ \ }}**@{},
     <0mm,0mm>*{};<-2mm,9mm>*{^{I_1}}**@{},
 <0mm,0mm>*{};<-8mm,-5mm>*{}**@{-},
 <0mm,0mm>*{};<-4.5mm,-5mm>*{}**@{-},
 <0mm,0mm>*{};<-1mm,-5mm>*{\ldots}**@{},
 <0mm,0mm>*{};<4.5mm,-5mm>*{}**@{-},
 <0mm,0mm>*{};<8mm,-5mm>*{}**@{-},
      <0mm,0mm>*{};<0mm,-7mm>*{\underbrace{\ \ \ \ \ \ \ \ \ \ \ \ \ \ \
      }}**@{},
      <0mm,0mm>*{};<0mm,-10.6mm>*{_{J_1}}**@{},
 <13mm,5mm>*{};<13mm,5mm>*{\bullet}**@{},
 <13mm,5mm>*{};<5mm,10mm>*{}**@{-},
 <13mm,5mm>*{};<8.5mm,10mm>*{}**@{-},
 <13mm,5mm>*{};<13mm,10mm>*{\ldots}**@{},
 <13mm,5mm>*{};<16.5mm,10mm>*{}**@{-},
 <13mm,5mm>*{};<20mm,10mm>*{}**@{-},
      <13mm,5mm>*{};<13mm,12mm>*{\overbrace{\ \ \ \ \ \ \ \ \ \ \ \ \ \ }}**@{},
      <13mm,5mm>*{};<13mm,14mm>*{^{I_2}}**@{},
 <13mm,5mm>*{};<8mm,0mm>*{}**@{-},
 <13mm,5mm>*{};<12mm,0mm>*{\ldots}**@{},
 <13mm,5mm>*{};<16.5mm,0mm>*{}**@{-},
 <13mm,5mm>*{};<20mm,0mm>*{}**@{-},
     <13mm,5mm>*{};<14.3mm,-2mm>*{\underbrace{\ \ \ \ \ \ \ \ \ \ \ }}**@{},
     <13mm,5mm>*{};<14.3mm,-4.5mm>*{_{J_2}}**@{},
 \end{xy}
$$
where $\sigma(I_1\sqcup I_2)$ is the sign of the shuffle
$I_1\sqcup I_2=(1,\ldots, m)$.}

\bip

\Proof
Claim (i) follows from the fact that $\LieBi^!(m,n)=1_m\ot sgn_n[1-n]$
and
Remark~2.2. Claim 2 is a straightforward though tedious
graph translation of the initial term,
$$
(\LieBi^!)^\vee(m,n)
\stackrel{d}{\lon}
{\bigoplus}_{(m,n)-{\mathrm trees}\, T \atop |edge(T)|=1}
(\LieBi^!)^*\ot {\mathrm Det}(T),
$$
of the definition 2.1 of the differential $d$ in ${\mathbf D}\LieBi^!$.
\hfill $\Box$

\bip

\sip

{\bf 3.3. A geometric model for $\LieBi_\infty$ structures}.
Let $V$ be a finite dimensional graded vector space. Then the graded formal manifold,
$\cM$, modeled on the infinitesimal neighbourhood of $0$ in the vector space
$V\oplus V^*[1]$ has an odd symplectic form $\om$ induced from the natural
pairing $V\ot V^*[1]\rar k[1]$. In particular, the graded structure sheaf $\f_{\cM}$
on $\cM$ has a degree $-1$ Poisson bracket, $\{\ \bullet\ \}$, such that
$$
\{f\bullet g\}=(-1)^{|f||g|+|f|+|g|}\{g\bullet f\}
$$
and the Jacobi identity is satisfied. The odd symplectic manifold $(\cM, \om)$
has two particular Lagrangian submanifolds, $\caL \subset \cM$ and
$\Pi\caL \subset \cM$
associated with, respectively, the subspaces
$0\oplus V^*[1]\subset V\oplus V^*[1]$ and $V\oplus 0\subset V\oplus V^*[1]$.

\bip

{\bf 3.3.1 Proposition.} {\em A $\LieBi_\infty$ algebra structure in a graded vector space
$V$ is the same as a degree two smooth function $\Ga\in \f_{\cM}$ vanishing on $\caL\cup
\Pi \caL$ and
 satisfying the equation $\{\Ga\bullet \Ga\}=0$.}

 \sip

 \Proof The manifold $\cM$ is isomorphic to the total space of the shifted 
cotangent bundle,
 $T^*_M[1]$, of the manifold $M$ of Proposition~1.5.1. Hence smooth functions on $\cM$
 are the same as smooth
 polyvector fields on $M$, and the Poisson
 bracket $\{\ \bullet\ \}$ on $\cM$ is the same as the Schouten bracket on $M$.
 \hfill $\Box$

 \bip

{\bf 3.4. $\LieBi_\infty$ morphisms.}
Let $(V, \{\mu_{m,n}\})$ and $(V', \{\mu_{m,n}'\})$ be two
 $\LieBi_\infty$ algebras.

 \bip

 {\bf 3.4.1. Definition.}
 A $\LieBi_\infty$ {\em morphism} $F: V \rar V'$ is, by definition, a
 symplectomorphism,
$F: (\cM, \om) \rar (\cM', \om')$ such that $F(\caL)\subset \caL'$,  
$F(\Pi\caL)\subset \Pi\caL'$ and
 $F^*\Gamma'=\Gamma$.

\bip

Thus a  $\LieBi_\infty$ {morphism} $F: V \rar V'$ is a pair of collections of
linear maps,
$$
\hspace{7mm}
\left\{ F_{m,n}: \odot^m V \ot \wedge^n V^* \rar V'[-n]\right\}_{m\geq 1, n\geq 0}, \ \
\left\{ \bar{F}_{m,n}:
\odot^m V \ot \wedge^n V^* \rar {V'}^*[1-n]\right\}_{n\geq 0, m\geq 1}
\hspace{8mm}\ (\star)
$$
satisfying the system equations, $F^*(\om')=\om$ and $F^*\Ga'=\Ga$. In particular,
the equation   $F^*(\om')=\om$ says that the linear maps,
$$
F_{1,0}: (V,\mu_{1,1}) \rar (V',\mu'_{1,1}) \ \ \mbox{and}\ \
\bar{F}_{0,1}: (V^*,\mu^*_{1,1}) \rar ({V'}^*,{\mu'}^*_{1,1}),
$$
are morphisms of complexes, while the equation  $F^*(\om')=\om$ says that
the composition,
$$
 F_{0,1}^*\circ \bar{F}_{1,0}: V\lon  V
$$
is the identity map.

\bip

 {\bf 3.4.2. Definition.} A  $\LieBi_\infty$-{morphism} $F: V \rar V'$
is called a {\em quasi-isomorphism}\,  if the morphisms of complexes
$$
F_{1,0}: (V,\mu_{1,1}) \rar (V',\mu'_{1,1}) \ \ \mbox{and}\ \
\bar{F}_{0,1}: (V^*,\mu^*_{1,1}) \rar ({V'}^*,{\mu'}^*_{1,1})
$$
induce isomorphisms in cohomology.
\bip

{\bf 3.4.3. Remark.} One might get an impression that the 
notions introduced above make sense only for
{\em finite dimensional}\, $\LieBi_\infty$ algebras. However, this is no more than
an artifact of the geometric intuition we tried to rely on in our definitions.
In fact, everything above (and below) make sense for infinite dimensional 
representations as well. For example, one can replace $(\star)$ by 
$$
\left\{F_{m,n}: \odot^m V  \rar \wedge^n V\ot  V'[-n]\right\}_{m\geq 1, n\geq 0}, \ \
\left\{ \bar{F}_{m,n}:
\odot^m V\ot {V'}[n-1]   \rar  \wedge^n V\right\}_{n\geq 0, m\geq 1},
$$
and reinterpret the  equations defining $\LieBi_\infty$ morphism
accordingly. 
For example, with this reinterpretation
 it is the morphism $F_{1,0}\circ \bar{F}_{0,1}$ which is the 
identity map.

 \bip

 {\bf 3.4.4. Contractible and minimal  $\LieBi_\infty$-structures.} Let
$V$ be a graded vector space and $(\cM=V\oplus V^*[1],\om)$ the associated odd 
symplectic manifold (as in the Sect.\ 3.3).
There is a one-to-one correspondence between differentials, $d:V\rar V$, 
and quadratic degree 2 function, $\Gamma_{quad}$, on $(\cM, \om)$, 
vanishing on $\caL\cup \Pi\caL$ and satisfying $[\Gamma_{quad}\bullet 
\Gamma_{quad}]=0$.
If $H^*(V,d)=0$, the associated data, $(\cM, \om, \Ga_{quad})$,
 is called a {\em contractible  $\LieBi_\infty$-structure}\, on $V$.

\sip

A $\LieBi_\infty$ structure, $(\cM, \om, \Ga)$, on $V$ is called {\em minimal}\, if
$\Ga=0 \bmod I^3$, where $I$ is the ideal of the distinguished point, $\caL\cap
\Pi\caL$, in $\cM$. Put another way, the formal power series $\Ga$ in some
(and hence any) coordinate system on $\cM$ begins with cubic terms at least.

\bip

 {\bf 3.4.5. Theorem (homotopy classification of $\LieBi_\infty$-structures,
cf.\  \cite{Ko,Ko2}).} {\em
 Each $\LieBi_\infty$ algebra is isomorphic to the tensor product
of a contractible   $\LieBi_\infty$ algebra and a minimal one.} 

\bip

\Proof 
Let $(\cM,\om, \Gamma)$ be the geometric equivalent of any given   
$\LieBi_\infty$ algebra. To prove the statement we have to construct coordinates,
$(x^a, y^a, z^\al, \psi_a, \phi_a, \xi_\al)$, on $\cM$ such that
\Bi
\item[(i)] $\om= 
\underset{\om_1}{\underbrace{(\sum_a(dx^a\wedge d\psi_\al + dy^a\wedge d\phi_\al)}}
 +\underset{\om_2}{\underbrace{\sum_\al dz^\al\wedge
d\xi_\al}}$,
\item[(ii)] $\caL$ is given by $x^a=y^a=z^\al=0$ while $\Pi\caL$ is given by
$\psi_a=\phi_a=\xi_\al=0$,
\item[(iii)] $\Gamma= \underset{\Ga_1}{\underbrace{\sum_a y^a\psi_a}} + 
\underset{\Ga_2}{\underbrace{\Phi(z^\al,\xi_\al)}}$ for some 
formal power series  $\Phi(z^\al,\xi_\al)$ which begins with cubic terms at least.
\Ei
For then 
$(\cM, \om, \Ga)\simeq (\cM_1, \om_1, \Ga_1)
\times  (\cM_2, \om_2, \Ga_2)$
with the first factor being a contractible $\LieBi_\infty$ structure
while the second factor a minimal one.

\sip

We shall establish existence of the above coordinates by induction. 

\sip

As the first step of the induction procedure we choose  arbitrary linear coordinates,
$\{t^A\}$, $A\in \{1, \ldots, \dim V\}$, on $V$ and the associated dual coordinates,
$\{\chi_A\}$, $|\chi_A|=-|t^A|+1$, on $V^*[1]$. The odd symplectic form
is given in these coordinates as $\om=\sum_A dt^A\wedge d\chi_A$,  
$\caL$ is given by $t^A=0$ while $\Pi\caL$ is given by
$\chi_A=0$. Then,
$$
\Gamma= \sum_{A,B} C_B^A t^A\chi_B \bmod I^3,
$$
for some constants $C_B^A$ which are nothing but the coefficients
of the differential, $d:V\rar V$, associated to the quadratic bit of $\Gamma$.
As we work over the field of characteristic zero, we can choose a cohomological
decomposition of $V$ with respect to this differential,
$$
V= H(V,d) \oplus B\oplus B[-1],
$$
so that $d$ vanishes on $H(V,d)$ and $B[-1]$ and, on the remaining summand $B$,
it  is equal to the natural
isomorphism $B\rar B[-1]$. Let $\{z^\al\}$ be some linear
coordinates in $H(V,d)$, $\{y^a\}$ linear coordinates on $B$ and $\{x^a\}$,
$|x^a|=|y^a|-1$, the associated (via the natural isomorphism) linear coordinates
in $B[-1]$. Denote by $(\xi_\al, \psi_a,\phi_a)$ the coordinates
on $V^*[1]$ dual to $(z^\al, x^a, y^a)$. In the resulting coordinate system
on $\cM$  the conditions (i)-(ii) are satisfied, while the condition 
(iii) is satisfied modulo  $I^3$.

\sip

Assume by induction that we have constructed a coordinate system on $\cM$ in which
conditions (i)-(iii) are satisfied $\bmod I^{N+1}$. Then we have,
$$
 \Gamma= \underset{\Ga_1}{\underbrace{\sum_a y^a\psi_a}}\ \ \  + 
\underset{{\rm polynomials\ of\atop degrees\ from\ 3\ to \ N}}
{\underbrace{\Phi^{\leq N}(z^\al,\xi_\al)}} \ \ 
+ \underset{{\rm polynomial\ of\ degree}\ N+1}
{\underbrace{\Gamma^{N+1}(x,y,z,\psi,\phi,\chi)}} \mod I^{N+2}.
$$
The equation $[\Ga\bullet \Ga]=0\bmod I^{N+2}$ implies,
$$
\delta\Gamma^{N+1}=0
$$
where $\delta$ is the following differential on $\f_\cM$,
$$
\Ba{rccc}
\delta: & \f_\cM & \lon & \f_\cM\\
&            f   & \lon & [{\sum_a y^a\psi_a}\bullet f]. 
\Ea
$$
Let $ B\in \f_\cM$ be an arbitrary polynomial of degree  $N+1$ and with $|B|=1$.
It gives rise to a symplectomorphism, $F: \cM\rar\cM$, given as $\exp v_B$
with the vector field $v_B$  defined by  
$$
\Ba{rccc}
v_B: & \f_\cM & \lon & \f_\cM\\
&            f   & \lon & [B\bullet f]. 
\Ea
$$
One has,
$$
F^*\Gamma =  {\sum_a y^a\psi_a} + {\Phi^{\leq N}(z^\al,\xi_\al)} +
\Gamma^{N+1} + \delta B \mod I^{N+2}.
$$
Thus $\Gamma^{N+1}(x,y,z,\psi,\phi,\chi)$ is a $\delta$-cycle defined up to
a $\delta$-coboundary. As cohomology of $\delta$ in $\f_\cM$ is equal to
$k[[z^\al, \xi_\al]]$, one can always find 
$\Gamma^{N+1}$ such that it is a function of $\{z^\al, \xi_\al\}$ only.
\hfill $\Box$

\bip

 {\bf 3.4.6. Corollary.} {\em If  $F: V \rar V'$ is a $Lie_1Bi_\infty$ quasi-isomorphism,
 then there exists a $Lie_1Bi_\infty$ quasi-isomorphism $G: V' \rar V$ such that
 on the cohomology level $[F_{1,0}]=[\bar{G}_{0,1}]^*$ and
  $[G_{1,0}]=[\bar{F}_{0,1}]^*$.
 }

\sip

{\bf Proof} is exactly the same as the proof of an analogous statement for $L_\infty$
algebras in \cite{Ko}.

\bip

{\bf 4. Minimal resolution of the operad ${\mathsf TF}$.}
By definition (see Sect.\ 1.4),  ${\mathsf TF}$ is a
quadratic dioperad
$$
{\mathsf TF}= \frac{ Free(E)}{Ideal <R>},
$$
where
\Bi
\item[(i)] $E(2,1):= k[\Sigma_2]\ot 1_1$ and  $E(1,2):= 1_1\ot 1_2[-1]$;
we represent
two basis vectors of $k[\Sigma_2]\ot 1_1$ by planar corollas
$$
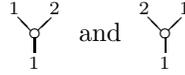
\begin{xy}
 <0mm,-0.55mm>*{};<0mm,-2.5mm>*{}**@{-},
 <0.5mm,0.5mm>*{};<2.2mm,2.2mm>*{}**@{-},
 <-0.48mm,0.48mm>*{};<-2.2mm,2.2mm>*{}**@{-},
 <0mm,0mm>*{\circ};<0mm,0mm>*{}**@{},
   <0mm,-0.55mm>*{};<0mm,-3.8mm>*{_1}**@{},
   <0.5mm,0.5mm>*{};<2.7mm,2.8mm>*{^2}**@{},
   <-0.48mm,0.48mm>*{};<-2.7mm,2.8mm>*{^1}**@{},
 \end{xy}
 \ \ \mbox{and} \ \
 \begin{xy}
 <0mm,-0.55mm>*{};<0mm,-2.5mm>*{}**@{-},
 <0.5mm,0.5mm>*{};<2.2mm,2.2mm>*{}**@{-},
 <-0.48mm,0.48mm>*{};<-2.2mm,2.2mm>*{}**@{-},
 <0mm,0mm>*{\circ};<0mm,0mm>*{}**@{},
   <0mm,-0.55mm>*{};<0mm,-3.8mm>*{_1}**@{},
   <0.5mm,0.5mm>*{};<2.7mm,2.8mm>*{^1}**@{},
   <-0.48mm,0.48mm>*{};<-2.7mm,2.8mm>*{^2}**@{},
 \end{xy}
 $$
and a basis vector of $E(1,2)$ by the symmetric corolla,
 $$
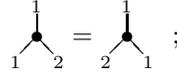
\begin{xy}
 <0mm,0.66mm>*{};<0mm,3mm>*{}**@{-},
 <0.39mm,-0.39mm>*{};<2.2mm,-2.2mm>*{}**@{-},
 <-0.35mm,-0.35mm>*{};<-2.2mm,-2.2mm>*{}**@{-},
 <0mm,0mm>*{\bullet};<0mm,0mm>*{}**@{},
   <0mm,0.66mm>*{};<0mm,3.4mm>*{^1}**@{},
   <0.39mm,-0.39mm>*{};<2.9mm,-4mm>*{^2}**@{},
   <-0.35mm,-0.35mm>*{};<-2.8mm,-4mm>*{^1}**@{},
\end{xy}
 =
 \begin{xy}
 <0mm,0.66mm>*{};<0mm,3mm>*{}**@{-},
 <0.39mm,-0.39mm>*{};<2.2mm,-2.2mm>*{}**@{-},
 <-0.35mm,-0.35mm>*{};<-2.2mm,-2.2mm>*{}**@{-},
 <0mm,0mm>*{\bullet};<0mm,0mm>*{}**@{},
   <0mm,0.66mm>*{};<0mm,3.4mm>*{^1}**@{},
   <0.39mm,-0.39mm>*{};<2.9mm,-4mm>*{^1}**@{},
   <-0.35mm,-0.35mm>*{};<-2.8mm,-4mm>*{^2}**@{},
\end{xy}\ \  ;
$$
\item[(ii)] the relations $R$ are generated by the following elements,
\Beqrn
 \begin{xy}
 <0mm,0mm>*{\bullet};<0mm,0mm>*{}**@{},
 <0mm,0.69mm>*{};<0mm,3.0mm>*{}**@{-},
 <0.39mm,-0.39mm>*{};<2.4mm,-2.4mm>*{}**@{-},
 <-0.35mm,-0.35mm>*{};<-1.9mm,-1.9mm>*{}**@{-},
 <-2.4mm,-2.4mm>*{\bullet};<-2.4mm,-2.4mm>*{}**@{},
 <-2.0mm,-2.8mm>*{};<0mm,-4.9mm>*{}**@{-},
 <-2.8mm,-2.9mm>*{};<-4.7mm,-4.9mm>*{}**@{-},
    <0.39mm,-0.39mm>*{};<3.3mm,-4.0mm>*{^3}**@{},
    <-2.0mm,-2.8mm>*{};<0.5mm,-6.7mm>*{^2}**@{},
    <-2.8mm,-2.9mm>*{};<-5.2mm,-6.7mm>*{^1}**@{},
 \end{xy}
\ + \
 \begin{xy}
 <0mm,0mm>*{\bullet};<0mm,0mm>*{}**@{},
 <0mm,0.69mm>*{};<0mm,3.0mm>*{}**@{-},
 <0.39mm,-0.39mm>*{};<2.4mm,-2.4mm>*{}**@{-},
 <-0.35mm,-0.35mm>*{};<-1.9mm,-1.9mm>*{}**@{-},
 <-2.4mm,-2.4mm>*{\bullet};<-2.4mm,-2.4mm>*{}**@{},
 <-2.0mm,-2.8mm>*{};<0mm,-4.9mm>*{}**@{-},
 <-2.8mm,-2.9mm>*{};<-4.7mm,-4.9mm>*{}**@{-},
    <0.39mm,-0.39mm>*{};<3.3mm,-4.0mm>*{^2}**@{},
    <-2.0mm,-2.8mm>*{};<0.5mm,-6.7mm>*{^1}**@{},
    <-2.8mm,-2.9mm>*{};<-5.2mm,-6.7mm>*{^3}**@{},
 \end{xy}
\ + \
 \begin{xy}
 <0mm,0mm>*{\bullet};<0mm,0mm>*{}**@{},
 <0mm,0.69mm>*{};<0mm,3.0mm>*{}**@{-},
 <0.39mm,-0.39mm>*{};<2.4mm,-2.4mm>*{}**@{-},
 <-0.35mm,-0.35mm>*{};<-1.9mm,-1.9mm>*{}**@{-},
 <-2.4mm,-2.4mm>*{\bullet};<-2.4mm,-2.4mm>*{}**@{},
 <-2.0mm,-2.8mm>*{};<0mm,-4.9mm>*{}**@{-},
 <-2.8mm,-2.9mm>*{};<-4.7mm,-4.9mm>*{}**@{-},
    <0.39mm,-0.39mm>*{};<3.3mm,-4.0mm>*{^1}**@{},
    <-2.0mm,-2.8mm>*{};<0.5mm,-6.7mm>*{^3}**@{},
    <-2.8mm,-2.9mm>*{};<-5.2mm,-6.7mm>*{^2}**@{},
 \end{xy}
\  \ \ &\in & Free(E)(1,3)  \\
&&\\
 \begin{xy}
 <0mm,2.47mm>*{};<0mm,-0.5mm>*{}**@{-},
 <0.5mm,3.5mm>*{};<2.2mm,5.2mm>*{}**@{-},
 <-0.48mm,3.48mm>*{};<-2.2mm,5.2mm>*{}**@{-},
 <0mm,3mm>*{\circ};<0mm,3mm>*{}**@{},
  <0mm,-0.8mm>*{\bullet};<0mm,-0.8mm>*{}**@{},
<0mm,-0.8mm>*{};<-2.2mm,-3.5mm>*{}**@{-},
 <0mm,-0.8mm>*{};<2.2mm,-3.5mm>*{}**@{-},
     <0.5mm,3.5mm>*{};<2.8mm,5.7mm>*{^2}**@{},
     <-0.48mm,3.48mm>*{};<-2.8mm,5.7mm>*{^1}**@{},
   <0mm,-0.8mm>*{};<-2.7mm,-5.2mm>*{^1}**@{},
   <0mm,-0.8mm>*{};<2.7mm,-5.2mm>*{^2}**@{},
\end{xy}
\  - \
\begin{xy}
 <0mm,-1.3mm>*{};<0mm,-3.5mm>*{}**@{-},
 <0.38mm,-0.2mm>*{};<2.2mm,2.2mm>*{}**@{-},
 <-0.38mm,-0.2mm>*{};<-2.2mm,2.2mm>*{}**@{-},
<0mm,-0.8mm>*{\circ};<0mm,0.8mm>*{}**@{},
 <2.4mm,2.4mm>*{\bullet};<2.4mm,2.4mm>*{}**@{},
 <2.5mm,2.3mm>*{};<4.4mm,-0.8mm>*{}**@{-},
 <2.4mm,2.5mm>*{};<2.4mm,5.2mm>*{}**@{-},
     <0mm,-1.3mm>*{};<0mm,-5.3mm>*{^2}**@{},
     <2.5mm,2.3mm>*{};<5.1mm,-2.6mm>*{^1}**@{},
    <2.4mm,2.5mm>*{};<2.4mm,5.7mm>*{^2}**@{},
    <-0.38mm,-0.2mm>*{};<-2.8mm,2.5mm>*{^1}**@{},
    \end{xy}
\ - \
\begin{xy}
 <0mm,-1.3mm>*{};<0mm,-3.5mm>*{}**@{-},
 <0.38mm,-0.2mm>*{};<2.2mm,2.2mm>*{}**@{-},
 <-0.38mm,-0.2mm>*{};<-2.2mm,2.2mm>*{}**@{-},
<0mm,-0.8mm>*{\circ};<0mm,0.8mm>*{}**@{},
 <2.4mm,2.4mm>*{\bullet};<2.4mm,2.4mm>*{}**@{},
 <2.5mm,2.3mm>*{};<4.4mm,-0.8mm>*{}**@{-},
 <2.4mm,2.5mm>*{};<2.4mm,5.2mm>*{}**@{-},
     <0mm,-1.3mm>*{};<0mm,-5.3mm>*{^1}**@{},
     <2.5mm,2.3mm>*{};<5.1mm,-2.6mm>*{^2}**@{},
    <2.4mm,2.5mm>*{};<2.4mm,5.7mm>*{^2}**@{},
    <-0.38mm,-0.2mm>*{};<-2.8mm,2.5mm>*{^1}**@{},
    \end{xy}
\  - \
\begin{xy}
 <0mm,-1.3mm>*{};<0mm,-3.5mm>*{}**@{-},
 <0.38mm,-0.2mm>*{};<2.2mm,2.2mm>*{}**@{-},
 <-0.38mm,-0.2mm>*{};<-2.2mm,2.2mm>*{}**@{-},
<0mm,-0.8mm>*{\circ};<0mm,0.8mm>*{}**@{},
 <-2.4mm,2.4mm>*{\bullet};<-2.4mm,2.4mm>*{}**@{},
 <-2.5mm,2.3mm>*{};<-4.4mm,-0.8mm>*{}**@{-},
 <-2.4mm,2.5mm>*{};<-2.4mm,5.2mm>*{}**@{-},
     <0mm,-1.3mm>*{};<0mm,-5.3mm>*{^2}**@{},
     <-2.5mm,2.3mm>*{};<-5.1mm,-2.6mm>*{^1}**@{},
    <-2.4mm,2.5mm>*{};<-2.4mm,5.7mm>*{^1}**@{},
    <0.38mm,-0.2mm>*{};<2.8mm,2.5mm>*{^2}**@{},
    \end{xy}
    \  - \
\begin{xy}
 <0mm,-1.3mm>*{};<0mm,-3.5mm>*{}**@{-},
 <0.38mm,-0.2mm>*{};<2.2mm,2.2mm>*{}**@{-},
 <-0.38mm,-0.2mm>*{};<-2.2mm,2.2mm>*{}**@{-},
<0mm,-0.8mm>*{\circ};<0mm,0.8mm>*{}**@{},
 <-2.4mm,2.4mm>*{\bullet};<-2.4mm,2.4mm>*{}**@{},
 <-2.5mm,2.3mm>*{};<-4.4mm,-0.8mm>*{}**@{-},
 <-2.4mm,2.5mm>*{};<-2.4mm,5.2mm>*{}**@{-},
     <0mm,-1.3mm>*{};<0mm,-5.3mm>*{^1}**@{},
     <-2.5mm,2.3mm>*{};<-5.1mm,-2.6mm>*{^2}**@{},
    <-2.4mm,2.5mm>*{};<-2.4mm,5.7mm>*{^1}**@{},
    <0.38mm,-0.2mm>*{};<2.8mm,2.5mm>*{^2}**@{},
    \end{xy}
\ \ \ &\in& Free(E)(2,2).
\Eeqrn
\Ei
\bip

\sip

Proposition 1.4.1 follows immediately from the following

\bip
{\bf 4.1. Proposition.} {\em  The minimal resolution,
${\mathsf TF}_\infty$, of the dioperad
 ${\mathsf TF}$ can be described as follows:

\sip

(i) As a dioperad of graded vector spaces, ${\mathsf TF}_\infty=Free(E)$,
where the collection, $E=\{E(m,n)\}$,
of
$(\Sigma_m,\Sigma_n)$-modules is given by
$$
E(m,n):=\left\{
\Ba{cl}
k[\Sigma_2]\ot 1_n & {\mathrm if}\ m=2, n\geq 2; \\
1_n[-1]             &{\mathrm if} \ m=1, n\geq 2;\\
0             &{\mathrm otherwise}.
\Ea
\right.
$$
(ii) If we represent two basis elements of $E(2,n)$ by
planar $(2,n)$-corollas, and the basis element of $E(1,n)$ by planar $(1,n)$
corolla,
$$
 
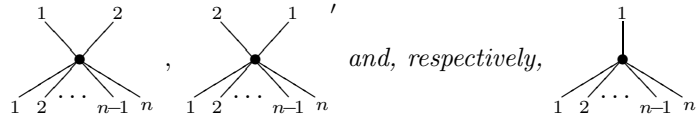
\begin{xy}
 <0mm,0mm>*{\bullet};<0mm,0mm>*{}**@{},
 <0mm,0mm>*{};<-4.5mm,5mm>*{}**@{-},
 <0mm,0mm>*{};<4.5mm,5mm>*{}**@{-},
   <0mm,0mm>*{};<-5mm,5.5mm>*{^1}**@{},
   <0mm,0mm>*{};<5.1mm,5.5mm>*{^2}**@{},
 <0mm,0mm>*{};<-8mm,-5mm>*{}**@{-},
 <0mm,0mm>*{};<-4.5mm,-5mm>*{}**@{-},
 <0mm,0mm>*{};<-1mm,-5mm>*{\ldots}**@{},
 <0mm,0mm>*{};<4.5mm,-5mm>*{}**@{-},
 <0mm,0mm>*{};<8mm,-5mm>*{}**@{-},
   <0mm,0mm>*{};<-8.5mm,-6.9mm>*{^1}**@{},
   <0mm,0mm>*{};<-5mm,-6.9mm>*{^2}**@{},
   <0mm,0mm>*{};<4.5mm,-6.9mm>*{^{n\hspace{-0.5mm}-\hspace{-0.5mm}1}}**@{},
   <0mm,0mm>*{};<9.0mm,-6.9mm>*{^n}**@{},
 \end{xy}
\ , \
 \begin{xy}
 <0mm,0mm>*{\bullet};<0mm,0mm>*{}**@{},
 <0mm,0mm>*{};<-4.5mm,5mm>*{}**@{-},
 <0mm,0mm>*{};<4.5mm,5mm>*{}**@{-},
   <0mm,0mm>*{};<-5mm,5.5mm>*{^2}**@{},
   <0mm,0mm>*{};<5.1mm,5.5mm>*{^1}**@{},
 <0mm,0mm>*{};<-8mm,-5mm>*{}**@{-},
 <0mm,0mm>*{};<-4.5mm,-5mm>*{}**@{-},
 <0mm,0mm>*{};<-1mm,-5mm>*{\ldots}**@{},
 <0mm,0mm>*{};<4.5mm,-5mm>*{}**@{-},
 <0mm,0mm>*{};<8mm,-5mm>*{}**@{-},
   <0mm,0mm>*{};<-8.5mm,-6.9mm>*{^1}**@{},
   <0mm,0mm>*{};<-5mm,-6.9mm>*{^2}**@{},
   <0mm,0mm>*{};<4.5mm,-6.9mm>*{^{n\hspace{-0.5mm}-\hspace{-0.5mm}1}}**@{},
   <0mm,0mm>*{};<9.0mm,-6.9mm>*{^n}**@{},
 \end{xy}
'
\ \mbox{and, respectively,}\
 \begin{xy}
 <0mm,0mm>*{\bullet};<0mm,0mm>*{}**@{},
 <0mm,0mm>*{};<0mm,5mm>*{}**@{-},
   <0mm,0mm>*{};<0mm,5.5mm>*{^1}**@{},
 <0mm,0mm>*{};<-8mm,-5mm>*{}**@{-},
 <0mm,0mm>*{};<-4.5mm,-5mm>*{}**@{-},
 <0mm,0mm>*{};<-1mm,-5mm>*{\ldots}**@{},
 <0mm,0mm>*{};<4.5mm,-5mm>*{}**@{-},
 <0mm,0mm>*{};<8mm,-5mm>*{}**@{-},
   <0mm,0mm>*{};<-8.5mm,-6.9mm>*{^1}**@{},
   <0mm,0mm>*{};<-5mm,-6.9mm>*{^2}**@{},
   <0mm,0mm>*{};<4.5mm,-6.9mm>*{^{n\hspace{-0.5mm}-\hspace{-0.5mm}1}}**@{},
   <0mm,0mm>*{};<9.0mm,-6.9mm>*{^n}**@{},
 \end{xy}
$$
with symmetric ingoing legs, then the differential
$d$ is given on generators by,
$$
d
\begin{xy}
 <0mm,0mm>*{\bullet};<0mm,0mm>*{}**@{},
 <0mm,0mm>*{};<0mm,5mm>*{}**@{-},
   <0mm,0mm>*{};<0mm,5.5mm>*{^1}**@{},
 <0mm,0mm>*{};<-8mm,-5mm>*{}**@{-},
 <0mm,0mm>*{};<-4.5mm,-5mm>*{}**@{-},
 <0mm,0mm>*{};<-1mm,-5mm>*{\ldots}**@{},
 <0mm,0mm>*{};<4.5mm,-5mm>*{}**@{-},
 <0mm,0mm>*{};<8mm,-5mm>*{}**@{-},
   <0mm,0mm>*{};<-8.5mm,-6.9mm>*{^1}**@{},
   <0mm,0mm>*{};<-5mm,-6.9mm>*{^2}**@{},
   <0mm,0mm>*{};<4.5mm,-6.9mm>*{^{n\hspace{-0.5mm}-\hspace{-0.5mm}1}}**@{},
   <0mm,0mm>*{};<9.0mm,-6.9mm>*{^n}**@{},
 \end{xy}
\ \ = \ \
 \sum_{{J_1\sqcup J_2=(1,\ldots,n)\atop
 {|J_1|\geq 2, |J_2|\geq 1}}
}\hspace{0mm}
 \begin{xy}
 <0mm,0mm>*{\bullet};<0mm,0mm>*{}**@{},
 <0mm,0mm>*{};<13mm,5mm>*{}**@{-},
 <0mm,0mm>*{};<-8mm,-5mm>*{}**@{-},
 <0mm,0mm>*{};<-4.5mm,-5mm>*{}**@{-},
 <0mm,0mm>*{};<-1mm,-5mm>*{\ldots}**@{},
 <0mm,0mm>*{};<4.5mm,-5mm>*{}**@{-},
 <0mm,0mm>*{};<8mm,-5mm>*{}**@{-},
      <0mm,0mm>*{};<0mm,-7mm>*{\underbrace{\ \ \ \ \ \ \ \ \ \ \ \ \ \ \
      }}**@{},
      <0mm,0mm>*{};<0mm,-10.6mm>*{_{J_1}}**@{},
 <13mm,5mm>*{};<13mm,5mm>*{\bullet}**@{},
 <13mm,5mm>*{};<13mm,10mm>*{}**@{-},
 <13mm,5mm>*{};<13mm,13mm>*{^1}**@{-},
 <13mm,5mm>*{};<8mm,0mm>*{}**@{-},
 <13mm,5mm>*{};<12mm,0mm>*{\ldots}**@{},
 <13mm,5mm>*{};<16.5mm,0mm>*{}**@{-},
 <13mm,5mm>*{};<20mm,0mm>*{}**@{-},
     <13mm,5mm>*{};<14.3mm,-2mm>*{\underbrace{\ \ \ \ \ \ \ \ \ \ \ }}**@{},
     <13mm,5mm>*{};<14.3mm,-4.5mm>*{_{J_2}}**@{},
 \end{xy}
$$
\Beqrn
d
\begin{xy}
 <0mm,0mm>*{\bullet};<0mm,0mm>*{}**@{},
 <0mm,0mm>*{};<-4.5mm,5mm>*{}**@{-},
 <0mm,0mm>*{};<4.5mm,5mm>*{}**@{-},
   <0mm,0mm>*{};<-5mm,5.5mm>*{^1}**@{},
   <0mm,0mm>*{};<5.1mm,5.5mm>*{^2}**@{},
 <0mm,0mm>*{};<-8mm,-5mm>*{}**@{-},
 <0mm,0mm>*{};<-4.5mm,-5mm>*{}**@{-},
 <0mm,0mm>*{};<-1mm,-5mm>*{\ldots}**@{},
 <0mm,0mm>*{};<4.5mm,-5mm>*{}**@{-},
 <0mm,0mm>*{};<8mm,-5mm>*{}**@{-},
   <0mm,0mm>*{};<-8.5mm,-6.9mm>*{^1}**@{},
   <0mm,0mm>*{};<-5mm,-6.9mm>*{^2}**@{},
   <0mm,0mm>*{};<4.5mm,-6.9mm>*{^{n\hspace{-0.5mm}-\hspace{-0.5mm}1}}**@{},
   <0mm,0mm>*{};<9.0mm,-6.9mm>*{^n}**@{},
 \end{xy}
\ \ & = & \ \
 \sum_{{J_1\sqcup J_2=(1,\ldots,n)\atop
 {|J_1|\geq 2, |J_2|\geq 0}}
}\hspace{0mm}
 \begin{xy}
 <0mm,0mm>*{\bullet};<0mm,0mm>*{}**@{},
 <0mm,0mm>*{};<13mm,5mm>*{}**@{-},
 <0mm,0mm>*{};<-8mm,-5mm>*{}**@{-},
 <0mm,0mm>*{};<-4.5mm,-5mm>*{}**@{-},
 <0mm,0mm>*{};<-1mm,-5mm>*{\ldots}**@{},
 <0mm,0mm>*{};<4.5mm,-5mm>*{}**@{-},
 <0mm,0mm>*{};<8mm,-5mm>*{}**@{-},
      <0mm,0mm>*{};<0mm,-7mm>*{\underbrace{\ \ \ \ \ \ \ \ \ \ \ \ \ \ \
      }}**@{},
      <0mm,0mm>*{};<0mm,-10.6mm>*{_{J_1}}**@{},
 <13mm,5mm>*{};<13mm,5mm>*{\bullet}**@{},
 <13mm,5mm>*{};<8.5mm,10mm>*{}**@{-},
 <13mm,5mm>*{};<17.5mm,10mm>*{}**@{-},
   <0mm,0mm>*{};<8.5mm,11mm>*{^1}**@{},
   <0mm,0mm>*{};<17.5mm,11mm>*{^2}**@{},
 <13mm,5mm>*{};<8mm,0mm>*{}**@{-},
 <13mm,5mm>*{};<12mm,0mm>*{\ldots}**@{},
 <13mm,5mm>*{};<16.5mm,0mm>*{}**@{-},
 <13mm,5mm>*{};<20mm,0mm>*{}**@{-},
     <13mm,5mm>*{};<14.3mm,-2mm>*{\underbrace{\ \ \ \ \ \ \ \ \ \ \ }}**@{},
     <13mm,5mm>*{};<14.3mm,-4.5mm>*{_{J_2}}**@{},
 \end{xy}
\\
&&
-\sum_{{J_1\sqcup J_2=(1,\ldots,n)\atop
 {|J_1|\geq 1, |J_2|\geq 1}}}\left(
 \begin{xy}
 <0mm,0mm>*{\bullet};<0mm,0mm>*{}**@{},
 <0mm,0mm>*{};<-4.5mm,5mm>*{}**@{-},
 <0mm,0mm>*{};<-4.5mm,6mm>*{^1}**@{},
 <0mm,0mm>*{};<13mm,5mm>*{}**@{-},
 <0mm,0mm>*{};<-8mm,-5mm>*{}**@{-},
 <0mm,0mm>*{};<-4.5mm,-5mm>*{}**@{-},
 <0mm,0mm>*{};<-1mm,-5mm>*{\ldots}**@{},
 <0mm,0mm>*{};<4.5mm,-5mm>*{}**@{-},
 <0mm,0mm>*{};<8mm,-5mm>*{}**@{-},
      <0mm,0mm>*{};<0mm,-7mm>*{\underbrace{\ \ \ \ \ \ \ \ \ \ \ \ \ \ \
      }}**@{},
      <0mm,0mm>*{};<0mm,-10.6mm>*{_{J_1}}**@{},
 <13mm,5mm>*{};<13mm,5mm>*{\bullet}**@{},
 <13mm,5mm>*{};<13mm,10mm>*{}**@{-},
  <13mm,5mm>*{};<13mm,11mm>*{^2}**@{},
 <13mm,5mm>*{};<8mm,0mm>*{}**@{-},
 <13mm,5mm>*{};<12mm,0mm>*{\ldots}**@{},
 <13mm,5mm>*{};<16.5mm,0mm>*{}**@{-},
 <13mm,5mm>*{};<20mm,0mm>*{}**@{-},
     <13mm,5mm>*{};<14.3mm,-2mm>*{\underbrace{\ \ \ \ \ \ \ \ \ \ \ }}**@{},
     <13mm,5mm>*{};<14.3mm,-4.5mm>*{_{J_2}}**@{},
 \end{xy}
 \ \ + \
  \begin{xy}
 <0mm,0mm>*{\bullet};<0mm,0mm>*{}**@{},
 <0mm,0mm>*{};<4.5mm,5mm>*{}**@{-},
 <0mm,0mm>*{};<4.5mm,6mm>*{^2}**@{},
 <0mm,0mm>*{};<-13mm,5mm>*{}**@{-},
 <0mm,0mm>*{};<-8mm,-5mm>*{}**@{-},
 <0mm,0mm>*{};<-4.5mm,-5mm>*{}**@{-},
 <0mm,0mm>*{};<-1mm,-5mm>*{\ldots}**@{},
 <0mm,0mm>*{};<4.5mm,-5mm>*{}**@{-},
 <0mm,0mm>*{};<8mm,-5mm>*{}**@{-},
      <0mm,0mm>*{};<0mm,-7mm>*{\underbrace{\ \ \ \ \ \ \ \ \ \ \ \ \ \ \
      }}**@{},
      <0mm,0mm>*{};<0mm,-10.6mm>*{_{J_1}}**@{},
 <-13mm,5mm>*{};<-13mm,5mm>*{\bullet}**@{},
 <-13mm,5mm>*{};<-13mm,10mm>*{}**@{-},
  <-13mm,5mm>*{};<-13mm,11mm>*{^1}**@{},
 <-13mm,5mm>*{};<-8mm,0mm>*{}**@{-},
 <-13mm,5mm>*{};<-12mm,0mm>*{\ldots}**@{},
 <-13mm,5mm>*{};<-16.5mm,0mm>*{}**@{-},
 <-13mm,5mm>*{};<-20mm,0mm>*{}**@{-},
     <-13mm,5mm>*{};<-14.3mm,-2mm>*{\underbrace{\ \ \ \ \ \ \ \ \ \ \ }}**@{},
     <-13mm,5mm>*{};<-14.3mm,-4.5mm>*{_{J_2}}**@{},
 \end{xy}
 \right)
\Eeqrn
}

\bip

\Proof Using criterion 2.4 it is easy to see that the dioperad ${\mathsf TF}$ is Koszul.
Then the cobar dual ${\mathbf D}{\mathsf TF}^!$ provides the required minimal resolution.
The rest is a straightforward calculation.
\hfill $\Box$

\bip

{\bf 5. A comment on ${\mathsf LieBi}_\infty$ algebras.} It was shown in \cite{G}
that the dioperad,
${\mathsf LieBi}$, of (usual) Lie bialgebras
 is Koszul so that its minimal resolution, ${\mathsf LieBi}_\infty$,
 can be constructed using the techniques
reviewed in Section 2. Here we  present its explicit graph description;
in fact, we prefer to show ${\mathsf LieBi}_\infty\langle 1\rangle$.

\bip

{\bf 5.1. Proposition.} {\em  The dioperad
 ${\mathsf LieBi}_\infty\langle 1\rangle$ can be described as follows.

(i) As a dioperad of graded vector spaces, ${\mathsf LieBi}_\infty
\langle 1 \rangle =Free(E)$,
where the collection, $E=\{E(m,n)\}$, of one dimensional
$(\Sigma_m,\Sigma_n)$-modules is given by
$$
E(m,n):=\left\{
\Ba{cl}
1_m\ot 1_n[2m-3] & {\mathrm if}\ m+n\geq 3; \\
0             &{\mathrm otherwise}.
\Ea
\right.
$$
(ii) If we represent a basis element of $E(m,n)$ by the unique
space $(m,n)$-corolla,
$$
 
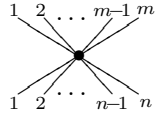
\begin{xy}
 <0mm,0mm>*{\bullet};<0mm,0mm>*{}**@{},
 <0mm,0mm>*{};<-8mm,5mm>*{}**@{-},
 <0mm,0mm>*{};<-4.5mm,5mm>*{}**@{-},
 <0mm,0mm>*{};<-1mm,5mm>*{\ldots}**@{},
 <0mm,0mm>*{};<4.5mm,5mm>*{}**@{-},
 <0mm,0mm>*{};<8mm,5mm>*{}**@{-},
   <0mm,0mm>*{};<-8.5mm,5.5mm>*{^1}**@{},
   <0mm,0mm>*{};<-5mm,5.5mm>*{^2}**@{},
   <0mm,0mm>*{};<4.5mm,5.5mm>*{^{m\hspace{-0.5mm}-\hspace{-0.5mm}1}}**@{},
   <0mm,0mm>*{};<9.0mm,5.5mm>*{^m}**@{},
 <0mm,0mm>*{};<-8mm,-5mm>*{}**@{-},
 <0mm,0mm>*{};<-4.5mm,-5mm>*{}**@{-},
 <0mm,0mm>*{};<-1mm,-5mm>*{\ldots}**@{},
 <0mm,0mm>*{};<4.5mm,-5mm>*{}**@{-},
 <0mm,0mm>*{};<8mm,-5mm>*{}**@{-},
   <0mm,0mm>*{};<-8.5mm,-6.9mm>*{^1}**@{},
   <0mm,0mm>*{};<-5mm,-6.9mm>*{^2}**@{},
   <0mm,0mm>*{};<4.5mm,-6.9mm>*{^{n\hspace{-0.5mm}-\hspace{-0.5mm}1}}**@{},
   <0mm,0mm>*{};<9.0mm,-6.9mm>*{^n}**@{},
 \end{xy}
$$
 then the differential
$d$ is given on generators by,
$$
d \begin{xy}
 <0mm,0mm>*{\bullet};<0mm,0mm>*{}**@{},
 <0mm,0mm>*{};<-8mm,5mm>*{}**@{-},
 <0mm,0mm>*{};<-4.5mm,5mm>*{}**@{-},
 <0mm,0mm>*{};<-1mm,5mm>*{\ldots}**@{},
 <0mm,0mm>*{};<4.5mm,5mm>*{}**@{-},
 <0mm,0mm>*{};<8mm,5mm>*{}**@{-},
   <0mm,0mm>*{};<-8.5mm,5.5mm>*{^1}**@{},
   <0mm,0mm>*{};<-5mm,5.5mm>*{^2}**@{},
   <0mm,0mm>*{};<4.5mm,5.5mm>*{^{m\hspace{-0.5mm}-\hspace{-0.5mm}1}}**@{},
   <0mm,0mm>*{};<9.0mm,5.5mm>*{^m}**@{},
 <0mm,0mm>*{};<-8mm,-5mm>*{}**@{-},
 <0mm,0mm>*{};<-4.5mm,-5mm>*{}**@{-},
 <0mm,0mm>*{};<-1mm,-5mm>*{\ldots}**@{},
 <0mm,0mm>*{};<4.5mm,-5mm>*{}**@{-},
 <0mm,0mm>*{};<8mm,-5mm>*{}**@{-},
   <0mm,0mm>*{};<-8.5mm,-6.9mm>*{^1}**@{},
   <0mm,0mm>*{};<-5mm,-6.9mm>*{^2}**@{},
   <0mm,0mm>*{};<4.5mm,-6.9mm>*{^{n\hspace{-0.5mm}-\hspace{-0.5mm}1}}**@{},
   <0mm,0mm>*{};<9.0mm,-6.9mm>*{^n}**@{},
 \end{xy}
\ \ = \ \
 \sum_{I_1\sqcup I_2=(1,\ldots,m)\atop {J_1\sqcup J_2=(1,\ldots,n)\atop
 {|I_1|\geq 0, |I_2|\geq 1 \atop
 |J_1|\geq 1, |J_2|\geq 0}}
}\hspace{0mm}
 \begin{xy}
 <0mm,0mm>*{\bullet};<0mm,0mm>*{}**@{},
 <0mm,0mm>*{};<-8mm,5mm>*{}**@{-},
 <0mm,0mm>*{};<-4.5mm,5mm>*{}**@{-},
 <0mm,0mm>*{};<0mm,5mm>*{\ldots}**@{},
 <0mm,0mm>*{};<4.5mm,5mm>*{}**@{-},
 <0mm,0mm>*{};<13mm,5mm>*{}**@{-},
     <0mm,0mm>*{};<-2mm,7mm>*{\overbrace{\ \ \ \ \ \ \ \ \ \ \ \ }}**@{},
     <0mm,0mm>*{};<-2mm,9mm>*{^{I_1}}**@{},
 <0mm,0mm>*{};<-8mm,-5mm>*{}**@{-},
 <0mm,0mm>*{};<-4.5mm,-5mm>*{}**@{-},
 <0mm,0mm>*{};<-1mm,-5mm>*{\ldots}**@{},
 <0mm,0mm>*{};<4.5mm,-5mm>*{}**@{-},
 <0mm,0mm>*{};<8mm,-5mm>*{}**@{-},
      <0mm,0mm>*{};<0mm,-7mm>*{\underbrace{\ \ \ \ \ \ \ \ \ \ \ \ \ \ \
      }}**@{},
      <0mm,0mm>*{};<0mm,-10.6mm>*{_{J_1}}**@{},
 <13mm,5mm>*{};<13mm,5mm>*{\bullet}**@{},
 <13mm,5mm>*{};<5mm,10mm>*{}**@{-},
 <13mm,5mm>*{};<8.5mm,10mm>*{}**@{-},
 <13mm,5mm>*{};<13mm,10mm>*{\ldots}**@{},
 <13mm,5mm>*{};<16.5mm,10mm>*{}**@{-},
 <13mm,5mm>*{};<20mm,10mm>*{}**@{-},
      <13mm,5mm>*{};<13mm,12mm>*{\overbrace{\ \ \ \ \ \ \ \ \ \ \ \ \ \ }}**@{},
      <13mm,5mm>*{};<13mm,14mm>*{^{I_2}}**@{},
 <13mm,5mm>*{};<8mm,0mm>*{}**@{-},
 <13mm,5mm>*{};<12mm,0mm>*{\ldots}**@{},
 <13mm,5mm>*{};<16.5mm,0mm>*{}**@{-},
 <13mm,5mm>*{};<20mm,0mm>*{}**@{-},
     <13mm,5mm>*{};<14.3mm,-2mm>*{\underbrace{\ \ \ \ \ \ \ \ \ \ \ }}**@{},
     <13mm,5mm>*{};<14.3mm,-4.5mm>*{_{J_2}}**@{},
 \end{xy}
$$

}

\bip

Let $V$ be a graded vector space, and let $\cM$ be the graded formal manifold
isomorphic to the neighbourhood of zero in $V[1] \oplus V^*[1]$. The manifold
$\cM$ has a natural even symplectic structure $\om$ induced from the paring
$V[1]\ot  V^*[1]\rar k[2]$; it also has two particular Lagrangian submanifolds,
 $\caL'$ and $\caL''$, modeled
on the subspaces $0\oplus  V^*[1] \subset V[1] \oplus V^*[1]$ and, respectively,
 $V\oplus 0 \subset V[1] \oplus V^*[1]$. The symplectic
form induces degree $-2$ Poisson bracket, $\{\ ,\ \}$, on the structure sheaf,
$\f_{\cM}$, of smooth functions on $\cM$.

\sip

The following result has been  independently obtained by Lyubashenko \cite{Lyu}.

\bip

{\bf 5.2 Corollary.} {\em A ${\mathsf LieBi}_\infty$ algebra 
structure in a graded vector space
$V$ is the same as a degree 3 smooth function $f\in \f_{\cM}$ vanishing on 
$\caL'\cup \caL''$ and
 satisfying the equation $\{f, f\}=0$.}

 \bip

${\mathsf LieBi}_\infty$ morphisms and quasi-isomorphisms are defined exactly as in 3.4.1
and 3.4.2; then an obvious  analogue of Theorem~3.4.5 holds true. We omit the details.

\bip

 \bip

  {\small

  \end{document}